\documentclass{amsart}

\usepackage{amsmath}
\usepackage{amssymb}
\usepackage{amscd}
\usepackage{amsthm}

 \numberwithin{equation}{section}

\theoremstyle{plain}

\newtheorem{Th}[equation]{Theorem}
\newtheorem{Prop}[equation]{Proposition}
\newtheorem{Le}[equation]{Lemma}
\newtheorem{Cor}[equation]{Corollary}

\theoremstyle{remark}

\theoremstyle{definition}

\newtheorem{Def}[equation]{Definition}

\newcommand{\ra}{\rightarrow}

\newcommand{\E}{\mathbb E}

\def\e{\emph}
\def\i{\infty}
\def\p{\partial}
\def\b{\begin}
\def\ra{\rightarrow}

\def\e{\emph}
\def\ra{\rightarrow}

\newcommand{\ol}{\overline}

\begin{document}

%\title{AMS Journal Sample}

\title{Tits Boundary of  $CAT(0)$ 2-complexes}
%\author{\flushleft{Xiangdong Xie}}

\author{Xiangdong Xie}
\address{Department of Mathematics, Washington University,
St.Louis, MO 63130.}
\email{xxie@math.wustl.edu}

\thanks{2000 \e{Mathematics Subject Classification.} Primary 20F67, 20F65; Secondary 57M20,  53C20.}
\thanks{\e{Key words and phrases.} {Tits boundary, Tits metric, $CAT(0)$, 2-complex, quasi-isometry, quasi-flat.}}

%\date{  }
%\maketitle

%\subjclass{Primary 20F67, 20F65; Secondary 57M20,53C20}

%\date{March 5, 2003.}

%\dedicatory{This paper is dedicated to our authors.}

%\keywords{Tits boundary, Tits metric, $CAT(0)$, 2-complex, quasi-isometry, quasi-flat.}

%\noindent
%{\bf{Mathematics Subject Classification(2000).}} Primary 53C23, 51F99, 57M20; 
% Secondary  53C70, 57M60.   
% \newline
%{\bf{Keywords.}} Tits boundary, Tits metric, CAT(0), 2-complex, geodesic.

%\keywords    test  \endkeywords
%\subjclass  51  \endsubjclass

\pagestyle{myheadings}

\markboth{{\upshape Xiangdong Xie}}{{\upshape Tits Boundary of $CAT(0)$ 2-complexes}}

%\footnotetext{\e{Mathematics Subject Classification(2000).} Primary 20F67, 20F65; Secondary 57M20,  53C20.}
      
%\keywords    test  \endkeywords
%\subjclass  51  \endsubjclass
\vspace{3mm}

\begin{abstract}{We investigate the Tits boundary   
of  locally compact  $CAT(0)$   $2$-complexes.    % By a result of B. Kleiner
    %small metric balls in    the Tits boundary   of $X$    are $R$-trees.
     In particular  we show    that  away from the endpoints,  
     a geodesic segment in the Tits boundary
  is the ideal boundary of an isometrically embedded Euclidean sector.
  % We   discuss   whether there is a lower bound  for   the 
 %distance between branch points   and  whether the length spectrum of
%simple closed geodesics in the Tits boundary is discrete.
   As  applications,  we provide   sufficient conditions for two points
 in the Tits boundary to   be the endpoints of a geodesic  in the  $2$-complex  and for a group
 generated by two hyperbolic isometries to contain a free group.  We also show  that  if two $CAT(0)$ 
  $2$-complexes are quasi-isometric then the cores  of their Tits boundaries are 
  bi-Lipschitz.}
\end{abstract}

% In this paper we investigate the Tits boundary of locally compact  CAT(0) 2-complexes
 % which admit cocompact group actions.   Small metric balls in the Tits boundary of such
%2-complexes are $R$-trees  by a result of B.Kleiner. We discuss whether   branch points 
%in the Tits boundary could accumulate and whether the length spectrum of
%simple closed geodesics in the Tits boundary is discrete.   
 %The starting point of this study is an   understanding of geodesic segments in
  %the Tits   boundary.}
%\tableofcontents  
 
\maketitle

\section{Introduction}

In this  paper   we    study  the Tits boundary of locally compact $CAT(0)$ 2-complexes. The $CAT(0)$ 2-complexes in this      paper     are  $CAT(0)$ piecewise Riemannian 2-complexes that admit cellular, isometric and cocompact   actions. Here each    closed  $2$-cell  is equipped with a Riemannian metric so that 
   it  is convex and its boundary is a broken geodesic. The metric on the 2-complex  is  the induced path metric.

Hadamard manifolds   are simply connected complete Riemannian manifolds of nonpositive sectional curvature. 
$CAT(0)$  spaces are counterparts of Hadamard manifolds in the category of  metric spaces.
%CAT(0)  spaces are metric versions of Hadamard manifolds, or simply connected complete Riemannian manifolds of %nonpositive sectional curvature.  
A $CAT(0)$ space is a complete simply
connected geodesic metric space so that all its  triangles   are at least as thin as the 
triangles in the Euclidean  space. 
$CAT(0)$  spaces  have many of the geometric properties enjoyed by Hadamard manifolds including 
 convexity of distance functions, uniqueness of geodesic segments and contractibility.  As in the case of Hadamard manifolds, a $CAT(0)$ space $X$ has a well-defined ideal boundary $\partial_\infty X$.  There is a topology on $\partial_\infty X$ called the cone topology and a metric   $d_T$  on $\partial_\infty X$ called the Tits metric.
 The topology induced by  $d_T$    is usually different from the cone topology.

 Given  a $CAT(0)$ space $X$, the Tits boundary $\p_T X$ means the ideal boundary  equipped 
with the Tits metric  $d_T$. 
 $\p_T X$   reflects the large scale geometry of $X$.  In particular,   $\p_T X$
  encodes information on 
%specifically 
large   flat subspaces of $X$  as well as  amount of negative curvature in $X$ (hyperbolicity).    
Tits  metric  and  Tits boundary    are  closely related to many 
 interesting questions in geometry.    They play an important role  in many  rigidity theorems
 (J. Heber \cite{H}, W. Ballmann \cite{B}, G. Mostow \cite{M}, B. Leeb \cite{L}).   
They are also very important  to the question  of whether a group
 of isometries of $X$ contains   a free 
group of rank two (K. Ruane \cite{R}, X. Xie \cite{X1}, \cite{X2}.)   

 %K. Ruane's theorem 
%provided a   sufficient condition in terms of the Tits metric for  a group generated 
 %by two hyperbolic isometries 
 %to  contain a free group of rank two.   

For an arbitrary $CAT(0)$ space  $X$,   $\p_T X$  is  complete and 
 CAT(1) (see   Section \ref{ideal} or   \cite{B}).   The Tits boundary of the product 
     of two $CAT(0)$ spaces is the    spherical   join of the Tits boundaries of the     two factors.    
When a   locally compact  $CAT(0)$ space  $X$ admits a cocompact isometric  action,
B. Kleiner (\cite{K}) showed the geometric dimension of  $\p_T X$   is 1 less than the 
maximal dimension of isometrically embedded Euclidean spaces  in $X$. Aside from these   general 
results   not much   is known about the Tits boundary of $CAT(0)$ spaces. Tits boundary is not even  well understood   for Hadamard manifolds
   or    2-dimensional (nonmanifold) $CAT(0)$ spaces.  
Tits boundary is well understood    only  for 
  a few  classes of $CAT(0)$ spaces:   $\p_T X$ is discrete  if  $X$ is   a  Gromov hyperbolic $CAT(0)$ space; 
$\p_T X$ is  a spherical building   if    $X$   is  a 
   higher rank symmetric space or Euclidean building;   $\p_T X$ contains interval components
      if  $X$ is 
 the universal cover of a nonpositively curved   graph manifold  
   (S. Buyalo-V. Schroeder \cite{BS}, C. Croke-B. Kleiner \cite{CK2}); 
   the Tits boundary has  also    been   studied
   for  real analytic    Hadamard $4$-manifolds
 admitting cocompact actions (C. Hummel-V. Schroeder  \cite{HS1}, \cite{HS2}) and 
the universal covers  of certain  torus  complexes (C. Croke-B. Kleiner  \cite{CK1}).   

%S. Buyalo,  V. Schroeder,C. Croke,   B. Kleiner,C. Hummel, V. Schroeder,U. Abresch,  V. Schroeder,

For the study of Tits boundary it  is fair to assume that the  $CAT(0)$ space admits a cocompact   isometric  action:
without the cocompactness condition, K. Kawamura and F. Ohtsuka (\cite{KO}) showed that any partition of the circle into points and open intervals can occur as  the Tits boundary of a    $CAT(0)$      piecewise Riemannian    metric on the plane. 
     
Let $X$ be    a   locally compact  $CAT(0)$  $2$-complex  admitting  a cocompact    action 
  by      cellular isometries.
   B. Kleiner's theorem  (Theorem C of \cite{K})  implies  that 
  $\partial_T X$   has geometric dimension at most 1.  
   It follows that  any closed metric ball with radius $r$ ($r< \pi/2$) in $\partial_T X$
   is an $R$-tree.  For an $R$-tree we can talk about  geodesic segments and  branch points.  Our first goal is to understand   geodesic segments and branch points in the Tits boundary.

A  \emph{sector}   is  a closed 
convex subset of  the  Euclidean plane $\mathbb{R}^2$  whose boundary is the union of 
 two  rays emanating from the origin.  
   We equip a   sector   with the
 induced metric.     The image of  an isometric embedding from a sector  into a   $CAT(0)$ space   $X$  is called a \emph{flat sector}  in $X$. 
We notice   if $S$ is a flat sector in a  $CAT(0)$ space   $X$, then   the  Tits boundary  
$\partial_T S$    of   $S$   is a closed interval  and isometrically embeds  into  
   the Tits boundary  $\partial_T X$    of   $X$.
The  following theorem  says that away from the endpoints, a   segment   in the Tits boundary 
   is  %spanned by an  isometrically embedded Euclidean  sector.
the   Tits boundary of a flat     sector.

\vspace{3mm}

\noindent
{\bf {Theorem \ref{t13}.}}
 \emph{Let $X$ be    a   locally compact  $CAT(0)$  $2$-complex  admitting  a cocompact action by
   cellular isometries,    and 
 %and $\partial_T X$  the Tits boundary of $X$.  Let 
$\gamma:[0,h]\rightarrow   \partial_T X$    a  geodesic in    $\partial_T X$
     with length 
$h\le \pi$. Then for any $\epsilon > 0$, there exists 
%an isometrically  embedded  Euclidean  sector $S$ in $X$ such that the    ideal  boundary    of 
   a  flat sector  $S$ in $X$   with       $\p_T S=\gamma([\epsilon, h-\epsilon])$.}

\vspace{3mm}

It follows %from  Theorem \ref{t13}  
that (see Proposition \ref{branchpo})   
 a branch point in  $\partial_T X$  is represented by a ray where two
   flat sectors branch   off.   
  Theorem  \ref{t13}   can not be improved to include the case 
 $\epsilon=0$: let $X$ be the universal cover of a torus complex considered in  \cite{CK1},
 then $\p_T X$ contains interval components;   some of these intervals  are not the Tits  boundaries 
  of any flat sectors.

Let $X$ be    as   in  Theorem \ref{t13}.
   Then small metric balls in $\p_T X$ are $R$-trees.  
An $R$-tree may have lots of branch points.   One may  ask if 
 there is  a 
constant $c=c(X)>0$ such that the distance between any two branch points in  $\p_T X$
%the Tits boundary 
is at least $c$. 
%B. Kleiner's theorem  also  implies circles   in $\p_T X$    are simple closed geodesics, 
  %in the Tits boundary are locally   geodesic, 
   %hence are rectifiable.  
Another  interesting question concerning  the Tits boundary is whether the lengths of circles in $\p_T X$
form a discrete set.  Note  B. Kleiner's  theorem  also   implies circles   in $\p_T X$  
 are simple closed geodesics, 
  %in the Tits boundary are locally   geodesic, 
   hence are rectifiable. 
%The answer to theses questions is positive when the interior
 %angles of $2$-cells are rational angles:

\vspace{3mm}

\noindent
 {\bf  {Theorem \ref{circlelength}}    and     {Theorem  \ref{t15}.}}
\emph{ Let $X$ be as in Theorem  \ref{t13}.
Suppose  the interior angles of all the   closed  $2$-cells  of $X$ are 
rational multiples of $\pi$. Then there is a positive integer $m$  such that:\newline
\e{(i)}   each topological circle in   $\p_T X$   has length  an integral  multiple of $\pi/m$;\newline
\e{(ii)} the distance between any two branch points 
in $\p_T X$    is   either infinite or   an integral  multiple of  ${\pi}/{m}$.   }

\vspace{3mm}

 Given  a  geodesic  $c: R\ra X$ in a $CAT(0)$ space $X$, 
  the two points in the ideal boundary determined by $c_{[0,\i)}$ and $c_{(-\i,0]}$    are  
% denoted by $c(+\i)$ and  $c(-\i)$, and 
called the \e{endpoints} of $c$.   
As    an  application of Theorem \ref{t13} we  discuss when  
two points  $\xi, \eta\in \p_\i X$ 
      are   the endpoints of a geodesic in $X$.
Recall a necessary condition is   $d_T(\xi, \eta)\ge \pi$,  
 and  a sufficient condition is     $d_T(\xi, \eta)> \pi$. 
 We provide a criterion  for $\xi$    and  $\eta$ 
to be the endpoints of a geodesic in $X$ when
    $d_T(\xi, \eta)=\pi$
and $X$  is a $CAT(0)$ 2-complex. 

A point in the ideal  boundary of a $CAT(0)$ space is called a \e{terminal point}
 if it does not lie in the interior of any Tits geodesic.

\vspace{3mm}

\noindent
{\bf{ Theorem \ref{pivisi}.}}
\e{Let $X$ be as in Theorem  \ref{t13}.
%Let   $X$ be a    $CAT(0)$  2-complex    that  
 %admits a  proper, cocompact   action by cellular isometries.
If $\xi, \eta\in \p_\i X$   are not terminal points and   $d_T(\xi, \eta)\ge \pi$,
   then there is a geodesic
 in $X$  with    $\xi$ and $\eta$  as endpoints.}

\vspace{3mm}

%The conclusion  of Theorem \ref{pivisi}  does not hold  if $X$ is not a $CAT(0)$ $2$-complex.  For instance, 
 %the universal covers of certain real analytic $4$-manifolds with nonpositive sectional curvature
 %are counterexamples (\cite{HS1}, \cite{AS}).

As further application we provide a  sufficient  condition  for     a   group generated by two 
 hyperbolic isometries to contain a free group of rank two.  Recall an isometry
 $g:X\ra X$ is a   \e{hyperbolic isometry}   if there is a geodesic $c:R\ra X$  and a positive number
 $l$ such that $g(c(t))=c(t+l)$ for all $t$. The geodesic $c$  is called an \e{axis}   of $g$.
      Denote the two endpoints of   the axis   $c$ by  $g(+\i)$ and  $g(-\i)$.   
 
\vspace{3mm}

\noindent 
{\bf{Theorem   \ref{th4.11}.}}
\e{Let $X$ be as in Theorem  \ref{t13}
%Let   $X$ be a  $CAT(0)$   $2$-complex    that  
% admits a  proper, cocompact   action by cellular isometries 
so that  each closed  $2$-cell is isometric to a convex polygon in the Euclidean plane   with 
all   its interior 
angles     rational    multiples  of   $\pi$, and 
   $g_1$,  $g_2$  two hyperbolic isometries 
 of $X$.   Suppose  $g_1(+\i)$, $g_1(-\i)$, $g_2(+\i)$, $g_2(-\i)$  are not terminal points   and 
  $d_T(\xi, \eta)\ge \pi$ for 
 any $\xi\in \{g_1(+\i), g_1(-\i)\}$,
   $\eta\in \{g_2(+\i), g_2(-\i)\}$. 
 Then  the group   generated by  $g_1$  and  $g_2$   contains  a free group 
 of rank two.}

\vspace{3mm}
% suggested  using support sets to study quasi-flats in 2-complexes.

  It  is well-known that   a quasi-isometry between  two  Gromov hyperbolic 
 spaces induces  a homeomorphism between their boundaries. This 
property does not hold for  $CAT(0)$  spaces.   C. Croke and B. Kleiner (\cite{CK1}) constructed 
  two quasi-isometric $CAT(0)$ $2$-complexes $X_1$ and $X_2$ 
 such  that  $\p_\i X_1$ and $\p_\i X_2$  are not homeomorphic   with respect to  the cone topology.
%a finite    2-complex  with two 
%different nonpositively curved metrics $d_1, d_2$ such that the   ideal boundaries of their universal coverings $X_1$ and $X_2$ 
%are not homeomorphic in the cone topology. 
It is also clear from their proof that    
 $\partial_T X_1$ and     $\partial_T X_2$   are not isometric. 
 %Croke  and Kleiner showed 
%that some components of the Tits boundaries are closed intervals. They further showed that a closed interval component in 
%$\partial _T{X_1}$ can correspond to a single point component in $\partial _T{X_2}$.  
 In general   it  is   unclear whether 
  the Tits boundaries of two quasi-isometric  $CAT(0)$  spaces    are  homeomorphic.  %Next result shows that, 
 %for  $CAT(0)$  2-complexes, the bi-Lipschitz  type of the ``main'' part of the Tits boundary is a quasi-isometry invariant.  

For any    $CAT(0)$  $2$-complex  $X$,    
set 
$ Core(\partial_T  X)=\cup c$
 where $c$   varies  over all the topological circles in $\partial_T X$.  Let $d_c$ be the induced path metric  of $d_T$  on 
 $Core(\partial_T X)$.   % Then $d_c\ge d_T$  always holds.      
   %Recall by  Corollary \ref{c1}  $\ol B(\xi,r)\subset \p_T X$ is  an $R$-tree 
%for any $\xi\in  \p_T X$ and any  $r$  with $0<r<\pi/2$.   It follows that   
%$d_c(\xi,\eta)<\i$  if  and only if  $\xi$, $\eta$ lie in the same path component of $Core(\p_T X)$
 %  and in this case there is  a minimal Tits geodesic  contained in  $Core(\p_T X)$   and 
%connecting   $\xi$ and  $\eta$. In particular,   
  % $d_c(\xi,\eta)=d_T(\xi,\eta)$ if $d_c(\xi,\eta)<\i$.   
  For  $\xi,\eta\in Core(\partial_T  X)$,   
 Kleiner's theorem implies  that $d_c(\xi,\eta)=d_T(\xi,\eta)$ if 
$\xi$, $\eta$ lie in the same path component of $Core(\p_T X)$, and $d_c(\xi,\eta)=\i$ 
otherwise.  

\vspace{3mm}

\noindent 
 {\bf{Theorem  \ref{t17}.}}
\e{Let $X_1$ and $X_2$ be two    $CAT(0)$  $2$-complexes. If $X_1$ and $X_2$  
are $(L,A)$ quasi-isometric, then $Core(\partial_T X_1)$ and $Core(\partial_T X_2)$  are  $ L^2$-bi-Lipschitz with respect to the metric $d_c$. }

\vspace{3mm}

The  paper  is organized as follows. In Section 2, we recall basic facts about 
$CAT(0)$ spaces   and  
  Tits boundary.    
 In Section  3 we use   support set to prove Theorem \ref{t13}. 
In Section 4 we give some applications of Theorem \ref{t13}; 
in this section we  first record  several   results concerning flat sectors and rays in $CAT(0)$  2-complexes, then we prove  Theorems   \ref{circlelength},           \ref{t15}    and \ref{pivisi}.  
In Section 5 we dicuss when a group generated by two hyperbolic isometries contains a free group
(Theorem \ref{th4.11}).    In Section  6 we  study Tits boundaries of quasi-isometric 
  $CAT(0)$ $2$-complexes  (Theorem  \ref{t17}).

\noindent
\textbf{Acknowledgment.} \textit{The author would like to thank 
  Bruce Kleiner, for numerous suggestions and discussions, without which this work would never be possible. In particular,  he  suggested  using support sets to study   $CAT(0)$   2-complexes.}

%the notion of support set is due to him.}

\section{Preliminaries }

The reader  is referred  to \cite{B}, \cite{BBr}, \cite{BH} and  \cite{K}  for more details on the material in this section.

\subsection{$CAT(\kappa)$ Spaces} \label {catk}

 Let $X$ be a metric space. For any $x\in X$ and any $r>0$,  ${B}(x,r)=\{x'\in X: d(x,x')< r\}$   and 
   $\ol{B}(x,r)=\{x'\in X: d(x,x')\le r\}$  are  respectively the open and 
 closed  metric balls   with center $x$ and radius $r$.
For any subset $A\subset X$ and any $\epsilon>0$, 
 the $\epsilon$-neighborhood of $A$ is  
$N_\epsilon(A)=\{x\in X: d(x, a)\le \epsilon\; \text{for some } a\in A\}$.
   For any two subsets  $A, B\subset X$,  the Hausdorff distance between $A$  and  $B$ is 
 $d_H(A, B)=\inf \{\epsilon: A\subset N_\epsilon(B), B\subset N_\epsilon(A)\}$; 
  $d_H(A, B)$  is defined to be $\i$ if there is no   $\epsilon>0$   with
  $A\subset N_\epsilon(B)$  and    $B\subset N_\epsilon(A)$.

The \e{Euclidean cone} over  a metric space $X$ is the
 metric space   $C(X)$   defined as follows.
 As a set $C(X)=X\times [0,\i)/{X\times \{0\}}$.  %The image of $(x,t)$ is denoted $tx$.
  We use $tx$  to denote the image of $(x,t)$.    We define 
$d(t_1x_1, t_2x_2)=\sqrt{t_1^2+t_2^2-2t_1t_2\cos(d(x_1,x_2))}$ if $d(x_1, x_2)\le \pi$,
   and $d(t_1x_1, t_2x_2)=t_1+t_2$ if $d(x_1, x_2)\ge \pi$.
%$d(t_1x_1, t_2x_2)=t_1+t_2$ if $d(x_1, x_2)\ge \pi$,  and 
%$d(t_1x_1, t_2x_2)=\sqrt{t_1^2+t_2^2-2t_1t_2\cos(d(x_1,x_2))}$ if $d(x_1, x_2)\le \pi$.
  The point $O=X\times \{0\}$ is called the cone point
  of $C(X)$.

   % The Euclidean  cone  $C_r(X)$  over $X$ with radius $r>0$ is the 
%closed metric ball $\ol{B}(O,r)\subset C(X)$ where $O=X\times \{0\}$ is the cone point
 % of $C(X)$.

Let $X$  be  a 
 metric space.   A \emph{geodesic}  in $X$   is a continuous map
    $c: I\rightarrow X$  from an interval  $I$ into $X$  such that, for any point $t\in I$, there exists a neighborhood $U$ of $t$
  with  $d(c(s_1),c(s_2))=|s_1-s_2|$ for   all  $s_1,s_2 \in U$.   
If
 the above equality holds for all $s_1, s_2 \in I$,
then we call $c $ a \emph{minimal geodesic}.
 The image of a geodesic shall also be called a geodesic.   When $I$ is a closed interval
$[a,b]$, we say  $c$ is a \e{geodesic segment} of length $b-a$ and  $c$ connects  $c(a)$  and  $c(b)$.  
A metric space $X$ is called a \emph{geodesic metric space} if for any  two points  $x,y\in X$ there is a 
   minimal geodesic segment connecting them.  

A \emph{ triangle}  in  a metric 
space $X$  is the union of   three geodesic segments $c_i: [a_i, b_i]\rightarrow X$ ($i=1, 2, 3$)
where  $c_1(b_1)=c_2(a_2)$, $c_2(b_2)=c_3(a_3)$   and  
$c_3(b_3)=c_1(a_1)$.    For any real number $\kappa$,   let $M^2_{\kappa}$ stand for the 2-dimensional simply connected complete Riemannian manifold with constant   sectional 
curvature $\kappa$, and   $D(\kappa)$  denote  the diameter of   $M^2_{\kappa}$ ($D(\kappa)=\infty $ if $\kappa\le 0$).    Given a triangle $\Delta=c_1\cup c_2\cup c_3$
in    $X$ where $c_i: [a_i, b_i]\rightarrow X$ ($i=1, 2, 3$),    a triangle $\Delta'$ 
in  $M^2_{\kappa}$     is a  \emph{comparison triangle}    for $\Delta$ if they have the same edge lengths, that is,  if  $\Delta'=c_1'\cup c_2'\cup c_3'$
  and  $c_i': [a_i, b_i]\rightarrow  M^2_{\kappa}  $  ($i=1, 2, 3$).   A point $x'\in \Delta'$ corresponds to a point $x\in  \Delta$ if there is some  $i$ and  some $t_i\in [a_i, b_i]$ with $x'=c_i'(t_i)$ and  $x=c_i(t_i)$.
We notice if the perimeter of a triangle $\Delta=c_1\cup c_2\cup c_3$
in $X$ is less than $2D(\kappa)$, that is, if $length(c_1)+length(c_2)+
length(c_3) < 2D(\kappa)$,  then   there is a unique comparison 
triangle (up to isometry) in $M^2_{\kappa}$   for $\Delta$. 

%\begin{itemize}
%\item \pkg{gen-j} for journal articles (uses \pkg{amsart})
%\item \pkg{gen-p} for proceedings volumes (uses \pkg{amsproc})
%\item \pkg{gen-m} for monographs (uses \pkg{amsbook})
%\end{itemize}

\begin{Def}\label{d1}
{Let  $\kappa\in R$.  A complete metric  space $X$ is called  a \emph {$CAT(\kappa)$ space} if\newline
{(i)} Every two points $x, y\in X$ with $d(x, y)< D(\kappa)$ are connected  by a 
    minimal  geodesic segment;\newline
{(ii)} For any triangle  $\Delta$ in $X$ with perimeter less than $2D(\kappa)$ and any two points $x,y\in \Delta$, the inequality  $d(x,y)\le d(x',y')$ holds,  where  $x'$ and $y'$ are the points on a comparison triangle for $\Delta$ corresponding to $x$ and $y$  respectively.}
%A complete metric space has \emph{curvature $\le \kappa$ } if each point has a 
%$CAT(\kappa)$ neighborhood.}  

\end{Def}

When $X$ is a $CAT(\kappa)$ space and $x, y\in X$ with $d(x, y)< D(\kappa)$,  the definition above
 implies there is a unique  minimal geodesic segment   $c_{xy}: [0, d(x,y)]\ra X$ with 
  $c_{xy}(0)=x$, $c_{xy}(d(x,y))=y$.
 We use  $xy$ to denote the 
image of    $c_{xy}$.

% When      $x$ and $y$ are as in Lemma \ref{uniquegeo} we denote the 
%image of the 
%unique geodesic segment from $x$ to $y$ by $xy$.  By abuse of  language,
%we sometimes also call $xy$ the geodesic segment from $x$ to $y$.

  A \emph{finite metric graph} is a finite graph where each edge  has a positive length. 
   We equip   the graph    with  the induced path 
metric.    Notice that a  finite metric graph is CAT(1) if and only if it has no simple  
 loop   with  length   strictly  less than  $2\pi$.  

A  geodesic  metric space $X$ is an  \emph{$R$-tree} if for any 
  triangle  $\Delta=c_1\cup c_2\cup c_3$  in $X$, $c_1$ is contained in the union 
  $c_2\cup c_3$.  

 %three points
  %$x, y,z \in X$ and any three geodesic segments $c_{xy}$, $c_{yz}$, $c_{zx}$ between them,
  %$c_{xy}$ is contained in the union  $c_{yz}\cup c_{zx}$.
    
\begin{Def}\label{branchpt}
{Let $X$ be an $R$-tree.  A point $p\in X$ is called a \emph{branch point} of $X$ if
  $X-\{p\}$ has at least three components.}
\end{Def}

%Another  class of $CAT(\kappa)$ spaces is provided by   piecewise constant curvature complexes. Let  $M^n_{\kappa}$    be the   
%n-dimensional simply connected complete Riemannian manifold with constant sectional curvature $\kappa$. A 
%  \emph{cell} in $M^n_{\kappa}$  
%is a compact  intersection of finitely many half spaces
% in $M^n_{\kappa}$. A piecewise constant curvature $\kappa$ complex $X$ is a cell 
%complex formed by gluing together cells in  $M^n_{\kappa}$ via isometries of their faces.  We assume $X$ is locally finite and there is% a uniform positive lower bound for the diameters of cells in $X$. In this case the induced  path metric on $X$ is complete. We will al%ways equip $X$ with 
%the induced path metric. Notice for each vertex $v\in X$, the link
% $Link(X,v)$ is a piecewise spherical complex, i.e., a piecewise  constant curvature 1 complex.

%The following result tells us when a piecewise constant curvature  $\kappa$ complex still has curvature bounded above by   $\kappa$.

%\begin{Th}(Alexandrov, Gromov\cite{G}, Ballmann\cite{GH})\label{t7}
%\emph{Suppose $X$ is a piecewise constant curvature $\kappa$ complex. Then $X$ has curvature upper bound $\kappa$ if and only %if for each vertex $v\in X$, the link
% $Link(X,v)$ is CAT(1).  }
%\end{Th}

\subsection{Space of Directions} \label{direc}

{A \emph{pseudo-metric} on a set $X$ is a function $d:  X\times X \rightarrow [0,   \infty)$ that is symmetric and satisfies the triangle inequality.}

If $(X,d)$ is a  pseudo-metric space, then we get a metric space $(X^*,  d^*)$ by letting $ X^*$ be the set of maximal zero
 diameter subsets and  setting  $ d^*(S_1,S_2):=d(s_1,s_2)$ for any $s_i\in S_i$.  

Let $X$ be a $CAT(\kappa)$ space.    If  $p,x,y\in X$  and $d(p,x)+d(x,y)+d(y,p)< 2 D(\kappa)$, then there is a well-defined geodesic triangle    $\bigtriangleup pxy$.   
                 The \emph{comparison angle}  of the triangle $\bigtriangleup pxy$ at $p$ is defined
to be the angle of the comparison triangle in $M^2_\kappa$ for $\bigtriangleup pxy$ at the vertex corresponding to $p$;  this angle is denoted   by  $\widetilde {\angle_p}(x,y)$.   The $CAT(\kappa)$  condition implies 
that if $x'\in {px}-\{p\}$ and  $y'\in {py}-\{p\}$  then 
    $\widetilde {\angle_p}(x',y')     \le \widetilde {\angle_p}(x,y)$.   Therefore if we let $x'\in {px}$,  $y'\in {py}$   tend to $p$
then   $\widetilde {\angle_p}(x',y') $  has a limit;   we call this   limit  the angle between ${px}$ and ${py}$   at $p$,  and denote it by $\angle_p(x,y)$.  If we let   $y'\in  {py}$   tend to $p$,    then   $\widetilde {\angle_p}(x,y') $  also tends to $\angle_p(x,y)$.  
The function  $p\rightarrow \angle_p(x,y)$   is upper semi-continuous.    $\angle_p$ defines a  pseudo-metric on the collection of geodesic segments leaving $p$.    We define $\Sigma^*_p X$ to be   the metric space
 associated  to the pseudo-metric  $\angle_p$. 
  The  space of directions at $p$ is the completion of  $\Sigma^*_p X$, and is denoted   by  $\Sigma_p X$.
For any $x\in B(p, D(\kappa))-\{p\}$, the point in  $\Sigma_p X$  coming from the geodesic segment $px$ shall be called the initial direction 
of $px$, and denoted  by $\log_{p}(x)$.  Thus we have a 
map $\log_p:        B(p, D(\kappa)) -\{p\}\rightarrow  \Sigma_p X$.

\begin{Th}\emph{(I. Nikolaev \cite{N})\label{t8}}
{Let $X$ be a $CAT(\kappa)$ space and $p\in X$.   Then $\Sigma_p X$ is a CAT(1) space.}
\end{Th}

\subsection{$CAT(0)$ $2$-Complexes} \label{2complex}

A $2$-dimensional $CW$-complex   is called a \e{polygonal complex}  if\newline
(1) all the attaching maps  are homeomorphisms;\newline
(2) the intersection of any two closed cells  is either empty or exactly
 one closed cell.

A $0$-cell  is also called a vertex.  
  %A   closed $2$-cell whose boundary contains $n$ vertices  is called a $n$-gon. 

    A polygonal complex  is  
 \emph{piecewise Riemannian }    if the following conditions hold:\newline
(1) For each closed 2-cell   $A$,  there   are    $n$ ($n\ge 3$) points $v_1, \cdots, v_n\in \p A$ 
and  a   Riemannian metric  on $A$ such that $v_iv_{i+1}$ ($i$ mod $n$) is a geodesic segment in the 
 Riemannian metric and the interior angle at $v_i$ ($1\le i\le n$) is at most $\pi$;\newline
%it is convex and its boundary 
%is a broken geodesic;\newline
(2) For any two closed 2-cells $A_1$ and $A_2$  with $A_1\cap A_2\not=\phi$, 
  the metrics on $A_1$ and $A_2$ agree
when restricted     to    $A_1\cap A_2$.

Let  $X$ be a piecewise Riemannian  polygonal complex and $x\in X$.   The link 
$Link(X, x)$  is a  metric graph  defined as follows.
    Let $A$ be a closed $2$-cell containing $x$. The  unit tangent space $S_x A$  of $A$ at $x$  is 
isometric to the unit circle with length $2\pi$.  We first define a subset 
  $Link(A,x)$  of $S_x A$. 
 For any $v\in S_x A$, 
  $v\in Link(A,x)$   if and only if 
the initial segment of the geodesic with initial point $x$ and initial direction $v$ lies in
   $A$.  Then $Link(A, x)=S_x A$  if $x$ lies in the interior of $A$;  $Link(A,x)$ is a 
  closed semicircle (with length $\pi$)
 if $x$ lies in the interior of a $1$-cell contained in $A$; and 
   $  Link(A,x)$  is a closed segment with length $\alpha$
 if $x$ is a vertex of $A$ and the interior angle 
 of $A$ at $x$ is $\alpha$.   Similarly  if $x$ is contained in a closed $1$-cell $B$ we can define 
  $S_x B$ and $Link(B, x)\subset S_x B$. We note $S_x B$ consists of two points at distance $\pi$
 apart, $Link(B, x)=S_x B$  if $x$ lies in    the  interior of $B$ and $Link(B, x)$ consists of a single point
 if $x$ is a vertex of $B$.  When $x$   lies  in a closed $1$-cell $B$ and $B$ is contained in a 
   closed $2$-cell $A$,  $S_x B$ and $Link(B, x)$ can be naturally identified with subsets of 
  $S_x A$ and $Link(A, x)$ respectively.

We define 
$Link(X, x)=\cup_A Link(A,x)$, where $A$ varies over all closed  $1$-cells and $2$-cells containing $x$. 
 Here $Link(B, x)$ is identified with a subset of $Link(A, x)$ as indicated in the last paragraph when
$x$   lies  in a closed $1$-cell $B$ and $B$ is contained in a 
   closed $2$-cell $A$.
  We let $d_x$ be   the induced path metric on $Link(X,x)$.

The following is a     corollary  of   Ballmann and Buyalo's Theorem (\cite{BBr}).   

\begin{Prop}   \label{cat0}
{Let $X$ be a   simply connected locally compact   piecewise Riemannian   polygonal complex  equipped with the induced path metric. 
 Suppose  $X$ admits a cocompact    action  by cellular isometries.   Then $X$ is a $CAT(0)$ space if and only if 
   the following conditions hold:\newline
\e{(i)}   the Gauss curvature of the open 2-cells is bounded from above by 0;\newline
\e{(ii)}  for every   vertex   $v$ of $X$  every simple loop in $Link(X, v)$ has length at least $2\pi$.}
\end{Prop}

A $CAT(0)$ $2$-complex in this paper shall always mean 
a polygonal complex satisfying all the conditions in Proposition \ref{cat0}.

Let $X$ be a $CAT(0)$ $2$-complex. 
 We notice for any $x\in X$, there is a natural identification between $\Sigma^*_x X=\Sigma_x X$ 
  and   $Link(X,x)$, and the path metric on $\Sigma_x X$ corresponds to $d_x$ on $Link(X,x)$.   

Let $X$ be a $CAT(0)$ space.   A  subset  $A\subset X$   is   a \e{convex} subset if  $xy\subset A$ 
 for any    $x, y\in A$.  Let $A\subset X$ be   a
 closed convex subset. 
 The orthogonal projection  onto $A$,  $\pi_A: X\ra A$  can be defined as follows:
  for any $x\in X$ the inequality $d(x, \pi_A(x))\le d(x, a)$  holds  for all 
  $a\in A$.  It follows that for any $x\notin A$ and any $a\not=\pi_A(x)$ we have
  $\angle_{\pi_A(x)}(x, a)\ge \pi/2$.   $\pi_A$  is 1-Lipschitz:
 $d(\pi_A(x), \pi_A(y))\le d(x,y)$
   for any $x,y\in X$.

Let  $X$  be  a  $CAT(0)$ space.   Then
  $\angle_x(y,z)+\angle_y(x,z)+\angle_z(x,y)\le \pi$
  for any three distinct points
$x,y,z\in X$.

\subsection{Ideal Boundary  of a  $CAT(0)$ Space} \label{ideal}

%We are mainly concerned with CAT(0) spaces. A basic feature of CAT(0) spaces is the convexity of distance function, that is,
%if $\sigma_1: [a,b]\rightarrow X$ and  $\sigma_2: [c,d]\rightarrow X$ are two geodesics in a CAT(0) space $X$ then the function
%$(s,t)\rightarrow d(\sigma_1(s),\sigma_2(t))$ is a convex function on $[a,b]\times [c,d]$.
%Thus given a CAT(0) space $X$, for any two points $x,y\in X$, there is a unique geodesic connecting 
%$x$ and $y$, and this geodesic varies continuously with respect to the two endpoints. 
%It follows that all the CAT(0) spaces are contractible.

% and all the spaces of nonpositive curvature are $K(\pi,1)$ spaces thanks to the following theorem:

%\begin{Th}\emph{(Alexandrov, Gromov\cite{G}, Ballmann\cite{GH})\label{t9}}
%{Let $\kappa\le 0$. If $Y$ is a space with curvature bounded above by  $\kappa$, then its universal cover is a $CAT(\kappa)$ space.}
%\end{Th}

Let $X$ be a $CAT(0)$ space.    A  \e{geodesic ray} in $X$  is a geodesic $c:[0,\infty) \ra X$.
Consider the set of geodesic rays in $X$. Two geodesic rays $c_1$ and $c_2$  are said to be \e{asymptotic}  if $f(t):=d(c_1(t)$, $c_2(t))$ is a bounded function.    It  is   easy to check that this defines an equivalence  relation. The set of equivalence classes is denoted  by  $\partial_{\infty}X$ and called the 
 \emph{ideal boundary} of $X$.  If $\xi\in \partial_{\infty}X$ and 
  $c$ is a geodesic ray  belonging to  $\xi$, we write $c(\infty)=\xi$.
For any $\xi\in \partial_{\infty}X$ and any $x\in X$, there is a unique geodesic ray 
$c_{x\xi}:[0,\infty)\rightarrow X$   with
$ c_{x\xi}(0)=x$ and $c_{x\xi}(\infty)=\xi$. 
The image of $c_{x\xi}$ is denoted  by  $x\xi$.

Set $\overline{X}=X\cup \partial_{\infty}X$.  The cone  topology  on 
$\overline{X}$  has   as a basis the open sets of $X$ together 
with the sets 
$$U(x,\xi,R,\epsilon)=\{z\in \overline{X} |z\notin B(x,R), 
d(c_{xz}(R),c_{x\xi}(R))<\epsilon\},$$
where $x\in X$, $\xi\in \partial_{\infty}X$ and $R>0$, $\epsilon>0$.  
The topology   on  $X$ induced   by   the cone topology  coincides with the 
metric topology on $X$.

%For example, the ideal boundary of a n-dimensional Hadamard manifold is homeomorphic to a (n-1)-sphere, and that of a regular  simplicial  tree with  valence $\ge 3$ is homeomorphic to a Cantor set.

We can also define a metric on  $\partial_{\infty}X$.
 Let  $c_1,c_2: [0, \infty) \rightarrow X$ be two geodesic rays 
    with $c_1(0)=c_2(0)=x$. For $t_1,t_2 \in (0,\infty)$,  consider the comparison angle    $\widetilde{\angle_x}(c_1(t_1),c_2(t_2))$. 
 The   $CAT(0)$ condition  implies that if $t_1 < t_1'$, $t_2 < t_2'$, then  $\widetilde{\angle_x}(c_1(t_1),c_2(t_2))\le 
\widetilde{\angle_x}(c_1(t_1'),c_2(t_2'))$. It follows that both     $\lim_{t\rightarrow 0}\widetilde{\angle_x}(c_1(t),c_2(t))$
 and $\lim_{t\rightarrow \infty}\widetilde{\angle_x}(c_1(t),c_2(t))$
 exist.  It can be proved that the limit 
$\lim_{t\rightarrow \infty}\widetilde{\angle_x}(c_1(t),c_2(t))$
 depends only on the points 
$\xi_1, \xi_2\in \p_\i X$ represented respectively
  by $c_1$ and $c_2$. We call 
  $\angle_T(\xi_1,\xi_2):=\lim_{t\rightarrow \infty}\widetilde{\angle_x}(c_1(t),c_2(t))$
 the \emph{Tits angle}    between $\xi_1$ and $\xi_2$.  The Tits metric $d_T$ on  
$\partial_{\infty}X$ is the path   metric induced by $\angle_T$. 
     We   denote
$\partial _T X:=(\partial_{\infty}X, d_T)$.
We  shall also call $\angle_x(\xi_1, \xi_2):=\lim_{t\rightarrow 0}\widetilde{\angle_x}(c_1(t),c_2(t))$
 the angle   at $x$   between $\xi_1$ and $\xi_2$.    Notice   $\angle_x(\xi_1, \xi_2)=\angle_x(p,q)$ 
for any $p\in x\xi_1$,  
$q\in x\xi_2$, $p,q\not= x$,  where $ \angle_x(p,q)$ is defined in   Section \ref{direc}. 
From the definition we see 
$$\angle_x(\xi_1, \xi_2)\le \widetilde{\angle_x}(c_1(t_1),c_2(t_2))\le \angle_T(\xi_1, \xi_2)\le d_T(\xi_1,\xi_2)$$
for all $t_1, t_2\in (0, \infty)$.

 It should be noted that the topology induced by $d_T$ is in general different from the cone topology. For instance, when $X=\mathbb{H}^n$, the n-dimensional real hyperbolic space, $\partial _T X$ is discrete while $\partial_{\infty}X$ with the cone topology is homeomorphic to $\mathbb{S}^{n-1}$. 
  
Here we record some basic properties of the Tits metric (see \cite{B} or \cite{BH}). 
For any geodesic $c: R\rightarrow X$ in a $CAT(0)$ space, we   call  the two points 
in $\partial_{\infty} X$ determined by the two rays $c_{|{[0, +\infty)}}$ and  $c_{|{(-\infty,0]}}$
 the endpoints of $c$, and denote them by $c(+\infty)$ and $c(-\infty)$  respectively.

\begin{Prop}\label{p2}
 {Let $X$ be a    locally   compact  $CAT(0)$   space, and $\xi_1$,   $\xi_2\in \partial_T X$. \newline
\e{(i)}  $\partial_T X$  is a CAT(1) space;\newline
\e{(ii)}  $\angle_T(\xi_1,\xi_2)=\sup_{x\in X} \angle_x(\xi_1,\xi_2)$;\newline
\e{(iii)}  If $d_T(\xi_1,\xi_2)> \pi$,  then there is a geodesic 
in $X$ with $\xi_1$ and $\xi_2$  as endpoints;\newline
\e{(iv)}  If   $d_T(\xi_1,\xi_2)< \i$,  then there is a   minimal  geodesic segment in $\partial _T X$ connecting $\xi_1$ and $\xi_2$.}  

\end{Prop}

   Recall a   \emph{sector}   is  a closed 
convex subset of  the  Euclidean plane $\mathbb{R}^2$   whose  boundary   is the union of 
 two  rays emanating from the origin.  
   We equip a   sector   with the
 induced metric.     The image of  an isometric embedding from a sector  into a   $CAT(0)$ space   $X$  is called a \emph{flat sector}  in $X$,  and the image of the origin is called the cone point of the flat sector.
We notice   if $S$ is a flat sector in a  $CAT(0)$ space   $X$, then 
$\partial_T S$ is a closed interval  and isometrically embeds  into  the Tits boundary  $\partial_T X$.

Let $X$ be a $CAT(0)$  space   and $p\in X$.  For each $\xi\in \partial_\i X$ the 
  geodesic ray 
$p\xi$ gives rise to a point in $\Sigma_p X$.
   Thus    $\log_p: X-\{p\}\rightarrow  \Sigma_p X$ extends to a    map 
$\overline X-\{p\}\rightarrow \Sigma_p X$, which is  continuous  and 
 shall still be denoted   by  $\log_p$.

\begin{Prop}\e{(\cite{B})}\label{p3}
{Let $X$ be a $CAT(0)$   space  and  $p\in X$. Then 
the  map $\log_p$ restricted to ${\partial_T X}$ is a 1-Lipschitz map,  that is, 
$\angle_p(\xi_1,\xi_2)\le d_T(\xi_1,\xi_2)$
  for  any $\xi_1, \xi_2$   in  $\p_T X$. % we have  $\angle_p(\xi_1,\xi_2)\le d_T(\xi_1,\xi_2)$. 
%the  map $\log_p$ restricted to ${\partial_T X}$ is a 1-Lipschitz map. 
If  $\xi_1, \xi_2   \in  \partial_T X$   with     
$d_T(\xi_1,\xi_2)=\angle_p(\xi_1,\xi_2)<\pi$, then the two rays $p\xi_1$, $p\xi_2$ bound a flat sector of angle $d_T(\xi_1,\xi_2)$.}  

\end{Prop}

{Let $Y$ be a $CAT(\kappa)$ space. We say the \emph{geometric dimension} of $Y$
 is $\le 1$ if $\Sigma_p Y$ is either empty or discrete for all $p\in Y$. The geometric dimension of $Y$ is defined to be 1 if it is $\le 1$ and $\Sigma_p Y$ is 
   nonempty   for   at least  one $p\in Y$.}
The following result is a special case of B. Kleiner's theorem (Theorem C of \cite{K}). The reader is referred to \cite{K} for the general result and general definition  of geometric dimension.  

\begin{Th}\emph{(B. Kleiner\cite{K})\label{t11}}
{Let $X$ be a    $CAT(0)$ 2-complex.% as defined in Section \ref{2complex}.
       Then the   geometric dimension of $  \partial_T X  $ is 
$\le 1$. 
If  furthermore  there exists an   isometric  embedding from the Euclidean plane  into $X$, then the geometric dimension of $  \partial_T X  $ is 1. }
\end{Th}

Let $X$ be a $CAT(0)$ 2-complex   %   as defined in Section \ref{2complex}  
     and   $a,b,c\in \partial_T X$   three distinct points  %in  $\partial_T X  $ 
     such  that  $d_T(a,b)$, $d_T(a, c) < {\pi}$.   Theorem \ref{t11}
     implies if $\log_a(b)\not=\log_a(c)$
  then  $\angle_a(b,c)=\pi$.  It follows that $ba\cup ac$ is  a geodesic segment
  when $\log_a(b)\not=\log_a(c)$.
  From this it   is   easy to derive the following corollary.

\begin{Cor}\label{c1}
{Let $X$ be a   $CAT(0)$ 2-complex.  %   as defined in Section \ref{2complex}. 
   Then for any  $\xi\in \partial_T X  $ and any $r: 0 < r < {\pi}/{2}$,
      the  closed metric ball 
   $\ol{B}(\xi, r)$ is an $R$-tree.}
\end{Cor}
%\begin{proof}     For any  $a,b,c \in  \partial_T X  $ with 
%$d_T(a,b)$, $d_T(a, c) < {\pi}$,    the intersection $ab\cap ac$ is a 
%geodesic segment (possibly  with length 0)   since  $\partial_T X  $  is a 
%CAT(1) space.  If  $d_T(a,b), d_T(a, c) \le  \frac{\pi}{2}$    
%and    $\log_a(b)\not=\log_a(c)$,
%then     $\angle_a(b,c)=\pi$ by Theorem \ref{t11}.  
%It follows that $ba\cup ac$ is still a geodesic segment. 
%Now              it  is   easy to see that for any  $\xi\in \partial_T X  $,  
%    any $r: 0 < r < \frac{\pi}{2}$ and any 
%three points $\xi_1, \xi_2, \xi_3\in  \bar B(\xi, r)$,  the 
%union    $\xi_1\xi_2\cup \xi_2\xi_3\cup \xi_3\xi_1$ is a geodesic 
%segment or a tripod.  It follows from  Lemma \ref{tree-tripod}  
% that     $\bar B(\xi, r)$ is    an $R$-tree.   

%\end{proof}

It follows from Corollary \ref{c1}   that each embedded path in  $\partial_T X  $   % the Tits boundary 
   is rectifiable and after reparameterization is  a geodesic.

\begin{Cor}  \label{r1}
{Let $X$ be a    $CAT(0)$ 2-complex.    %  as defined in Section \ref{2complex}.
  Then embedded paths  in $\partial_T X  $  are geodesics. In particular, all the topological circles
  in $\partial_T X  $ 
  are simple closed geodesics.}
\end{Cor}

%It follows from Corollary \ref{c1}   that each embedded path in the Tits boundary is rectifiable and after reparameterization they are locally distance minimizing. From now on we will assume all embedded paths   in the Tits boundary of a 2-complex  are parameterized by arc length.    In particular, all the topological circles in the Tits  boundary of a 2-complex are locally distance minimizing.}

\subsection{Quasi-isometry and Quasi-flats}\label{quaiso}

%A quasi-isometry embedding   between two metric spaces is a  Lipschitz map  for points not too
%close to each other:

\begin{Def}\label{d4}
{Let $L\ge 1$, $A\ge 0$.
A (not necessarily continuous) map $f:X\rightarrow Y$ between two metric spaces is called a 
  \emph{$(L,A)$  quasi-isometric embedding} if   
  the following holds    for all $x_1,x_2\in X$:  
$$\frac{1}{L}d(x_1,x_2)-A\le d(f(x_1), f(x_2))\le L d(x_1,x_2)+A.$$
  $f$  is   a \emph{$(L,A)$ quasi-isometry} if in addition $d_H(Y, f(X))\le A$.
   We  call $f$    a  quasi-isometric embedding if it is a $(L,A)$  quasi-isometric embedding
  for some $L\ge 1$, $A\ge 0$.}

%we have $Y=N_A(f(X))$, that is, the $A$-neighborhood of $f(X)$ is $Y$.}
\end{Def}

Notice if  $f:X\ra Y$ is  a $(L,A)$ quasi-isometry,  then there is some $A'>0$ 
  and a   $(L,A')$ quasi-isometry  $g:Y  \ra  X$    with   $d(g(f(x)), x)\le A'$,
    $d(f(g(y)), y)\le A'$ for  all $x\in X$, $y\in Y$. Such  a map $g$ is called a 
  \e{quasi-inverse}  of $f$.

Let  $X$  be a $CAT(0)$ $2$-complex.  
%A  quasi-isometric embedding from the real line into a metric space is called a \emph{quasi-geodesic}, and  
 The  image  of  a  quasi-isometric embedding from the  Euclidean plane  $\mathbb{R}^2$ %($k\ge 2$) 
    into   $X$    is   called a    \e{quasi-flat}.

%a metric space    $X$   is called a 
%\emph{quasi-$k$-flat}   in  $X$.  
%We also abuse language and say the image of the quasi-isometric embedding is a quasi-$k$-flat.
 % When $X$ is a $CAT(0)$ $2$-complex, then a quasi-2-flat in $X$ is also called a quasi-flat. 

\section{Segments in the Tits Boundary}

The main goal of this section is to establish the following result, which 
 says that away from the endpoints, a   segment   in the Tits boundary is
   the  Tits boundary  of   a flat   sector. The definition of a flat sector is given in  Section \ref{ideal}.

%by an  isometrically embedded Euclidean  sector.

\begin{Th} \label{t13}   
{Let $X$ be  a   $CAT(0)$   $2$-complex  as defined in Section \ref{2complex},   and 
$\gamma:[0,h]\rightarrow
\partial_T X$   a  geodesic    segment in the Tits 
boundary of $X$ with length 
$h\le \pi$. Then for any $\epsilon > 0$, there exists 
a    flat sector $S$ in $X$    with 
$\partial_T S=\gamma([\epsilon, h-\epsilon])$.}
 
\end{Th}

 %We observe it suffices to prove Theorem \ref{t13} for $h<\pi$. 

\subsection{A Gauss-Bonnet Type Theorem }  \label{gauss}

%                                 Let $X$ be a locally compact, simply connected polygonal complex. Assume each 2-cell is equip%ped with a Riemannian metric so that it is convex and  its boundary is a broken geodesic. 
%We equip $X$ with the induced path metric. Notice that the link at each 0-cell  is a finite metric graph.  Since the boundary %of each 2-cell is a broken geodesic,   $X$ is a CAT(0) space if and only if the sectional curvature of each 2-cell is $\le 0$ %and the link of each 0-cell is a CAT(1) space. We recall a finite metric graph is a CAT(1) space if and only 
%if it has no nontrivial loop of length $<$ $2\pi$.

In the proof of Theorem \ref{t13}   we will use a Gauss-Bonnet type theorem for noncompact piecewise Riemannian surfaces, which relates the Tits distance    (at infinity) and the curvature inside the surface.  We 
  now describe this result.

%Let $X$ be a piecewise Riemannian    2-complex (see definition \ref{polygona}),   and $v\in X$ a 0-cell. We denote the link of $X$ at $v$ by 
%$Link(X,v)$.  Notice $Link(X,v)$   is a finite metric graph where for each edge
 % $e$ of $Link(X,v)$  its length is the interior angle at $v$ 
%of the  2-cell corresponding to $e$.  
%The induced   path metric  in $Link(X,v)$   is denoted  $d_l$. 
 %  Recalling the definition of $\Sigma_v X$, we see $Link(X,v)$  and $\Sigma_v X$
%are homeomorphic  and 
 %$\angle_v(x_1,x_2)=\min\{\pi,d_l(\phi^v(x_1), \phi^v(x_2))\}$  for any $x_1, x_2\not= v$.
 %Notice that $d_l(\omega_1, \omega_2)$
%might be  $>  \pi$   while $\angle_v$ is always $\le \pi$.
 %$d_l(\omega_1, \omega_2)=\infty$ if $\omega_1$ and $\omega_2$ are in different components of $Link(X,v)$.  

Let  $F$ be a  piecewise Riemannian polygonal complex   with  the following properties:\newline
(1) $F$ is a  $CAT(0)$ space  with the induced path metric;\newline
(2) $F$ is   homeomorphic to a sector in the Euclidean  plane, in particular, $F$ is a manifold with boundary;\newline
(3) the manifold boundary of $F$ is   $\p F=c_1\cup c_2$, where $c_1$  and  $c_2$  are 
 geodesic rays  with  $c_1\cap c_2=\{p\}$, $p\in F$.

%We notice $F$ looks like a sector in the plane. 
   For each point $x$ in the interior of $F$, $Link(F,x)$ is a topological circle, and 
for each $x\in \p F$, $Link(F,x)$ is homeomorphic to a closed interval. 
 Let  $L(x)$  be  the  length of $Link(F,x)$.  The  deficiency   $k(x)$  at $x$ 
is defined as follows: $k(x):=2\pi-L(x)$   if $x$   lies  in the interior  of $F$,
 and $k(x):=\pi-L(x)$  if $x\in \p F$.     For each 2-cell $A$ of $F$, we let
$C(A)$ be the total curvature of $A$.  
  Set $e(F)=\sum_{A} C(A) + \sum_x k(x)$, 
where $A$  varies over all 2-cells 
and $x\not=p$  varies over all the 0-cells of $F$ different from $p$.

Let $\xi_1, \xi_2\in \partial_T F$ be represented by $c_1$ and $c_2$ respectively. 

\begin{Th}\e{(K. Kawamura, F. Ohtsuka    \cite{KO})}\label{t12}
{Let  $F$   and $\xi_1$, $\xi_2$ be as above. If $d_T(\xi_1,\xi_2)<\pi $, then $e(F)$ is  finite  and $e(F)=\angle_p(\xi_1, \xi_2)-d_T(\xi_1,\xi_2)$.}
\end{Th}

\subsection{Reduction} \label{reduction}

In this section we  reduce the proof of   Theorem \ref{t13} to the injectivity of a 
certain map.

        Let $X$ and  $\gamma:[0,h]\rightarrow \partial_T X$ be as in  
     Theorem \ref{t13}. 
Recall for each $x\in X$, 
there is a map
$\log_x:\overline X-\{x\}\rightarrow Link(X,x)$ which sends  $\xi\in     \overline X-\{x\}      $ 
   to
 the  initial direction of $x\xi$ 
at $x$.  For  $\xi_1$, $ \xi_2 \in \partial_T X$ with
 $d_T(\xi_1,\xi_2)< \pi$, let $\xi_1\xi_2\subset \p_T X$ be the 
unique geodesic
segment from $\xi_1$ to $\xi_2$ and $C(\xi_1 \xi_2)$ the 
Euclidean cone over $\xi_1\xi_2$.
Notice that $C(\xi_1 \xi_2)$ is a 
sector    with   angle
$d_T(\xi_1,\xi_2)$. Let $O$ be the cone point of $C(\xi_1 \xi_2)$.   For any
$x\in X$, there is a map $\rho_x:C(\xi_1 \xi_2)\rightarrow X$, 
where for any
$\xi \in \xi_1\xi_2$ the ray $O\xi$ is 
isometrically mapped to
the ray $x\xi$.  We notice $\rho_x$ is a $1$-Lipschitz   map.

For each 2-cell $A$ of $X$, let $C(A)$ be the total  curvature of $A$.  Note $C(A)\le 0$ since $A$ has
   nonpositive sectional curvature. 
Set 
$$\epsilon_1=\min\{-C(A):  A \text{  is a 2-cell with   }  C(A)\not=0\}.$$  $\epsilon_1$ is defined to be $\infty$ if there is no 
       $2$-cell  $A$ with 
$ C(A)\not=0$.   For a  finite metric graph $G$ and    an  edge path $c$ in
$G$, we denote the length of $c$ by $l(c)$.   Since $X$ is a  
     $CAT(0)$   $2$-complex, for any
 $x\in X$ each simple loop in $Link(X,x)$ has length at least $2\pi$.  
Define 
$$\epsilon_2=\min\{l(c)-2\pi: c \mbox{ is a simple loop in  } Link(X,x) 
   \mbox{ with }  l(c)\not=2\pi,   \;x\in X  \}.  $$  
$\epsilon_2$ is defined to be $\infty$ if  
there is no simple loop $c$ in any $Link(X,x)$ with  $l(c)\not=2\pi$.
Since $X$ admits a cellular cocompact    isometric   action, $\epsilon_1$  and $\epsilon_2$
 are  well-defined  and greater than $0$,  although  they   may be  $\infty$.

\begin{Le}\label{l8}
{Let  $x_0\in X$  with  $\angle_{x_0}(\gamma(0),\gamma(h))   \ge d_T(\gamma(0),\gamma(h))-\epsilon/4$,
   where $ 0<\epsilon<\min\{\epsilon_1, \epsilon_2,h\}     $.    If  
%Let $\epsilon $ be   a number with    $ 0<\epsilon<\min\{\epsilon_1, \epsilon_2,h\}     $. 
%If for some   $x_0\in X$  with  $\angle_{x_0}(\gamma(0),\gamma(h))   \ge d_T(\gamma(0),\gamma(h))-\epsilon/4$ 
the map   
${\log_{x_0}}_{|{\gamma(\epsilon/2)\gamma(h-\epsilon/2)}}$ is  injective, 
and  the surface 
$F_0:=\rho_{x_0}(C( \gamma(\epsilon/2)\gamma(h-\epsilon/2)) )$  is convex in $X$,  
then there
is a point $p\in X$ such that  %the map 
$\rho_p:  C(\gamma(\epsilon)\gamma(h-\epsilon) )\rightarrow X$ is an isometric embedding.}
\end{Le}
 
\begin{proof}  The   injectivity of the map ${\log_{x_0}}_{|{\gamma(\epsilon/2)\gamma(h-\epsilon/2)}}$
 implies that the map 
$$\rho_{x_0}:  C(\gamma(\epsilon/2)\gamma(h-\epsilon/2)) \rightarrow X$$
 is a topological embedding.  So
$F_0=\rho_{x_0}(C(\gamma(\epsilon/2)\gamma(h-\epsilon/2)) )$  is   homeomorphic to a sector
  in the Euclidean plane.  The manifold boundary $\p F_0$
of $F_0$ is   the  union of two geodesic rays 
  $x_0\gamma(\epsilon/2)$  and  $x_0\gamma(h-\epsilon/2)$   with   
  $x_0\gamma(\epsilon/2)\cap x_0\gamma(h-\epsilon/2)=\{x_0\}$.  
      $F_0$  is clearly a closed subset  of  $X$.

  Since  $F_0$ is closed and convex,  it is a $CAT(0)$ space with the induced path metric.
     Thus   we can apply Theorem \ref{t12}  to the surface $F_0$:
$$e(F_0)=\angle_{x_0}(\gamma(\epsilon/2),\gamma(h-\epsilon/2))  -
  d_T(\gamma(\epsilon/2),\gamma(h-\epsilon/2)),$$
 where  $e(F_0)=\sum_{A}C(A) +\sum_{p\not=x_0} k(p)$. By the  assumption on  $x_0$,    we  see 
the difference between $\angle_{x_0}(\gamma(\epsilon/2),\gamma(h-\epsilon/2))$ and $d_T(\gamma(\epsilon/2),\gamma(h-\epsilon/2))$ is less than any
positive   $-C(A)$   for 2-cells  $A$  of $X$ contained in  $F_0$.  It 
  follows that  $C(A)=0$  and  so  $A$ is flat  for any 2-cell  $A$ 
 of $X$ contained in the interior 
of $F_0$.   Similarly the difference between $\angle_{x_0}(\gamma(\epsilon/2),\gamma(h-\epsilon/2))$ and $d_T(\gamma(\epsilon/2),\gamma(h-\epsilon/2))$ is less than any
positive   $l(c)-2\pi$ for  simple loops $c$ in $Link(X,x)$, $x\in X$. 
Therefore $ Link(F_0,x)$     has length $2\pi$  for any  $x$ in the interior of $F_0$.

Since the space $X$ admits a   cellular cocompact   isometric  action, the sizes 
of the cells are bounded. Therefore  there   exists  some constant $r>0$ such  that the 
surface $F_1:=F_0- N_{r}(\p F_0)$ is flat, where  $N_{r}(\p F_0)$ denotes the $r$-neighborhood
of $ \p F_0\subset X$.    Now   notice   for any   $p\in F_1$  with $d(p, \p F_0)$ sufficiently large,
     $\rho_p( C(\gamma(\epsilon)\gamma(h-\epsilon) ))\subset F_1$ and the lemma follows. 

\end{proof}

\subsection{Segments in the Tits Boundary}

\vspace{5mm}

Let $X$ and  $\gamma:[0,h]\rightarrow \partial_T X$ be as in  
     Theorem \ref{t13},
   and  $\epsilon_1$ and $\epsilon_2$  as defined in Section \ref{reduction}.  %We  may assume $h<\pi$. 
We observe it suffices to prove Theorem \ref{t13} for $h<\pi$. 
Set $\epsilon_0=\min \{\epsilon_1, \epsilon_2, h\}$   and  
   denote   $\xi=\gamma(0)$, $\eta=\gamma(h)$. 
For any positive number $\epsilon<\epsilon_0$,   choose a point $x_0$  in  some   open 2-cell 
   $A$  of  $X$ such that 
$$\angle_{x_0}(\xi, \eta)
            \ge d_T(\xi,\eta)-\epsilon/4.$$ 
Notice Theorem  \ref{t13} 
  follows from   Lemma \ref{l8}    and  the following   theorem:

\begin{Th}\label{t16}
{Let $\epsilon$ and $x_0$ be    as above. Then the map 
  ${\log_{x_0}}_{|{\gamma(\epsilon/2)\gamma(h- \epsilon/2)}}$   is injective   and 
the surface 
$F_0:=\rho_{x_0}(C( \gamma(\epsilon/2)\gamma(h-\epsilon/2)) )$  is convex in $X$.}
\end{Th}

We shall consider the loop
 $c:=x_0\xi\cup x_0\eta \cup \xi\eta       \subset \overline X=X\cup \partial_{\infty} X$
  and its  support set  in $X$. Notice 
$c$ is indeed a loop in the cone topology of      $\overline X$. 
 Theorem \ref{t13} was  originally proved using a different method. The use of support set is suggested by B. Kleiner which 
greatly shortened  the argument.

Notice for any $x\in X-c$,   $\log_x(c)$   is  a loop in the link   $Link(X, x)$. 

\begin{Def}\label{supp}
{The \emph{support set} $supp(c)$ of $c$ is the set of $x\in X-c$ such that  $\log_x(c)$  
   represents a  nontrivial    class   in 
$H_1(Link(X, x))$.}
\end{Def}

%\b{remark}\label{moreg}
%\e{B. Kleiner's   support set
%can be defined for general homology classes and general 
 %polyhedral complexes.  Here we  are trying to take a short route.}
%\end{remark}

We shall first show that $supp(c)$ is   topologically a surface.  % From now on we   fix a homeomorphism
 %$C: \mathbb{S}^1\rightarrow c$ and set $c(t)=C(t)$ for $t\in \mathbb{S}^1$.

\begin{Le}\label{s1}
Let $x\in X-c$. Then $x\in supp(c)$ if and only if  $\log_x(c)$  is homotopic to a simple loop of length $2\pi$ 
in $Link(X,x)$. 
\end{Le}

\begin{proof}
Since $Link(X,x)$  is a finite graph, any simple loop in $Link(X,x)$   represents   a 
  nontrivial    class in 
$H_1(Link(X, x))$.   Thus one direction is clear.  Next we assume $x\in supp(c)$. 
    Note   $\log_x(c)$ 
is homotopically nontrivial   as  it is homologically nontrivial.
Since 
$x\notin x_0\xi$, we have  %$\angle_{x_0}(y_1, y_2) < \pi$ for any
$d_{x}(\log_x(y_1),  \log_x(y_2)) < \pi$ for any 
$y_1, y_2\in x_0\xi$.   The fact that $Link(X,x)$  is a CAT(1)  space implies
   the path $\log_x(x_0\xi)\subset Link(X,x)$ is  relative homotopic to the geodesic segment
from $a_1:=\log_x(x_0)$  to $b_1:=\log_x(\xi)$.  Similarly the path $\log_x(x_0\eta)$  
  is  relative homotopic to the geodesic segment
from $a_1$ to $c_1:=\log_x(\eta)$ and  $\log_x(\xi\eta)$  
is  relative homotopic to the geodesic segment
from $b_1$  to $c_1$.   Thus  $\log_x(c)$   is homotopic to the loop 
$a_1b_1* b_1c_1* c_1a_1$.  By considering the two ideal triangles 
$\Delta(xx_0\xi)$  and $\Delta(xx_0\eta)$ 
 we have 
   $\angle_{x_0}(x, \xi)+d_{x}(a_1, b_1)= \angle_{x_0}(x, \xi)+\angle_x(x_0,\xi)\le \pi$
%$d_{x_0}(\log_{x_0}(x), \log_{x_0}(\xi))+d_{x}(a_1, b_1)\leq \pi$
and $\angle_{x_0}(x, \eta)+d_{x}(a_1, c_1)= \angle_{x_0}(x, \eta)+\angle_x(x_0,\eta)\le \pi$.
%$d_{x_0}(\log_{x_0}(x), \log_{x_0}(\eta))+d_{x}(a_1, c_1)\leq \pi$.  
On the other hand, 
$$d_T(\xi,\eta)-\epsilon/4 \le \angle_{x_0}(\xi,\eta)\leq 
\angle_{x_0}(\xi, x)+\angle_{x_0}(x, \eta).$$
%d_{x_0}(\log_{x_0}(x), \log_{x_0}(\xi))+ d_{x_0}(\log_{x_0}(x), \log_{x_0}(\eta)).$$ 
  It follows that

\begin{tabular}{lll}
  &  &$d_{x}(a_1, b_1)+d_{x}(b_1, c_1)+ d_{x}(c_1, a_1)$\\
  & $\leq$ & $d_{x}(a_1, b_1)+ d_T(\xi, \eta)+ d_{x}(c_1, a_1)$\\
   & $\le$ & $\epsilon/4 +  \angle_{x_0}(\xi,\eta)   +d_{x}(a_1, b_1)+  d_{x}(c_1, a_1)$ \\
& $\le$ & $\epsilon/4 +   \angle_{x_0}(\xi, x)+\angle_{x_0}(x, \eta) +d_{x}(a_1, b_1)+  d_{x}(c_1, a_1)$ \\
  & $\leq$ & $\epsilon/4+2\pi$.
\end{tabular}  \newline
 Now the loop 
$a_1b_1* b_1c_1* c_1a_1$ is homotopically nontrivial in a CAT(1) finite  metric graph 
 with length at most $\epsilon/4+2\pi$.  The choice of $\epsilon$ implies there is no 
simple loop in $Link(X,x)$ with length  at most $2\pi+\epsilon/4$ but 
strictly   greater than   $2\pi$.  Therefore  $a_1b_1* b_1c_1* c_1a_1$ and 
$\log_x(c)$  are homotopic to a simple loop of length $2\pi$.  

\end{proof}

For any topological space $Y$  and $r>0$, let $C_r(Y)=Y\times [0, r]/(Y \times \{0\})$
  be the cone over $Y$ with radius $r$.  
Since $c$ is   a  circle,  $C_r(c)$   is  homeomorphic to the  closed unit disk
  in the plane.    Let   $\p C_r(c)$   be the boundary circle of  $C_r(c)$.

%For any $x\in X$, let $U(x)$ be the union of all the closed 2-cells  containing  $x$.
 %  Then $U(x)$ is a neighborhood of $x$ in $X$.  Set $r(x)=d(x, \partial U(x))$. 
%Then   for any  $r$   with  $ 0< r < r(x)$,   $\ol{B}(x,r)$     
 %   is   homeomorphic to the   cone   $C_r(\partial \ol B(x,r))$.

For   any  $x\in X-c$      and   any  $r$ with $0 < r <  d(x, x_0\xi\cup x_0\eta)$,
%$0 < r < \min\{r(x), d(x, x_0\xi\cup x_0\eta)\}$,   
   we can define a map 
$f_{x, r}: C_r(c)  \ra \ol B(x,r)$
  by  letting   $f_{x,r}(z,s)$ ($z\in c$, $0\le s\le r$)  be 
 the point on the geodesic $xz$ at distance $s$ from $x$.
  We observe that $f_{x,r}$  represents a class  in the relative homology group
     $H_2(X, X-\{x\})$.  By using homotopy along geodesic segments we see
  for $r_1$, $r_2$ with $0< r_1, r_2< d(x, x_0\xi\cup x_0\eta)$,
   $f_{x,r_1}$   and $f_{x,r_2}$  represent the same class in $H_2(X, X-\{x\})$.

For any $x\in X$, let $U(x)$ be the union of all the closed 2-cells  containing  $x$.
   Then $U(x)$ is a neighborhood of $x$ in $X$.  Set $r(x)=d(x, \partial U(x))$. 
%Then   for any  $r$   with  $ 0< r < r(x)$,   $\ol{B}(x,r)$     
 %   is   homeomorphic to the   cone   $C_r(\partial \ol B(x,r))$.

\begin{Le}\label{s2}
Let  $x\in X-c$   and  $r$ with 
$0 < r <d(x, x_0\xi\cup x_0\eta)$.
% \min\{r(x), d(x, x_0\xi\cup x_0\eta)\}$. 
Then     $x\in supp(c)$ if and only   if  $0\not=[f_{x,r}]\in H_2(X, X-\{x\})$.
\end{Le}

\begin{proof}    We may assume 
$0 < r < \min\{r(x), d(x, x_0\xi\cup x_0\eta)\}$
 by the remark preceding the lemma.
By excision we have $H_2(X,X-\{x\})\cong H_2(\ol{B}(x,r), \ol{B}(x,r)-\{x\})$. 
  %The  choice of $r$ implies $\ol B(x,r)$ is a cone over $\partial \ol B(x,r)$ and 
%$ \ol B(x,r)-x$ is homeomorphic to $\partial \ol B(x,r)\times (0, r]$. 
 The exact sequence 
for the pair $(\ol{B}(x,r), \ol{B}(x,r)-\{x\})$ implies the boundary homomorphism
$H_2(\ol{B}(x,r), \ol{B}(x,r)-\{x\})\rightarrow H_1(\ol{B}(x,r)-\{x\})$  is an isomorphism. 
  Notice the choice of $r$ implies 
the map $\log_x$ restricted to $ \partial \ol{B}(x,r)$  is a homeomorphism
  from $ \partial \ol{B}(x,r)$  to $Link(X, x)$.  It follows that
  $H_1(\partial \ol{B}(x,r))\cong H_1(Link(X,x))$.  Also notice 
 $H_1(\ol{B}(x,r)-\{x\})\cong H_1(\partial \ol{B}(x,r))$.  The composition of these 
 isomorphisms   is  an isomorphism $g: H_2(X,X-\{x\})\ra H_1(Link(X,x))$. 
  Now it  is   not hard to check that     $g$    maps the class $[f_{x,r}]$  to the class
 $[\log_x(c)]\in H_1(Link(X,x))$.

 %The image of $[f_{x,r}]$ under this homomorphism is the class in 
%$H_1(\ol B(x,r)-\{x\})$  represented by the map $c\rightarrow \partial \ol B(x,r)$ where 
%$z\in c$ is sent to $f_{x,r}(z,r)$. Now observe that the choice of $r$ implies 
%the map $\log_x$ is a homeomorphism when restricted to $ \partial \ol B(x,r)$. 
 % Now the lemma follows since the map $\log_x: c\rightarrow Link(X,x)$ factors 
%through $\partial \ol B(x,r)$.   

\end{proof}

    Lemma \ref{s1}   implies if   $x\in supp(c)$   then 
$\log_x(c)$  is homotopic to a simple loop $c_x$ in  $Link(X,x)$. 
%Recall the  edges in $Link(X,x)$  are in 1-1 correspondence with the 2-cells of 
%$X$ that are incident to $x$. 
  For each $x\in supp(c)$, let $S(x)$ be the union of all closed 2-cells 
%that correspond to the edges in
  that give rise to the simple loop    $c_x$.  % Since $c_x$ is a simple loop, $S(x)$ is 
 %homeomorphic   to the   closed unit disk in the  Euclidean   plane. 
Then $S(x)$ is  homeomorphic   to   the  closed   unit disk in the  Euclidean   plane
   and  contains $x$ in its interior
such that $Link(S(x),x)=c_x$ has length $2\pi$.

%\b{Le}\label{surface}
%{ Let  $c\subset \p_T X$ be a  simple  closed geodesic  with $length(c)<4\pi$.  Then for each $x\in supp(c)$,
 % there is some $r>0$   with 
  %$supp(c)\cap \ol{B}(x,r)=S(x)\cap \ol{B}(x,r)$.}
%\end{Le}

\begin{Le}\label{s3}
 For  any    $x\in supp(c)$,
  there is     $r>0$     such that  the  following    holds:
  $supp(c)\cap \ol{B}(x,r)=S(x)\cap \ol{B}(x,r)$.
   In particular   $supp(c)$ is a   $2$-dimensional manifold, and for each
$x\in supp(c)$ the circle  $Link(supp(c), x)$  has length $2\pi$.
\end{Le}

\begin{proof}
  For $x\in supp(c)$, %Lemma \ref{s1} says the loop $\log_x(c)$ is homotopic to a 
%simple loop $\gamma_x\subset Link(X,x)$ with length $2\pi$.  Let
%$D$ be the union of those closed 2-cells that give rise to the simple loop
%$\gamma_x$.  Then $D$ is a closed disk that contains $x$ in its interior
%such that $Link(D,x)=\gamma_x$ has length $2\pi$. 
 choose $r$ with  $0 < r < \frac{1}{2} \min\{r(x), d(x, x_0\xi\cup x_0\eta)\}$.
%Then $ \ol{B}(x,r)$ is a cone over the finite graph $\partial \ol B(x,r)$. 
Set $K=\ol{B}(x, r/2)$.
Notice $f_{x,r}$ represents a class in $H_2(X, X-K)$.   
%By Lemma \ref{simpleloop}  $\log_x(c)$ is homotopic to a simple closed geodesic 
 %in $Link(X,x)$. 
It  follows from the definitions of 
$S(x)$ and $f_{x,r}$  that 
$[f_{x,r}]=[S(x)\cap \ol{B}(x,r)]\in H_2(X, X-K) $.

For any $y\in \ol{B}(x,r/4)$,   by using homotopy along geodesic segments
 we see  $f_{y, r}$ and $f_{x, r}$  as maps  from $(C_r(c), \p C_r(c))$ 
  to 
 $(X,X-K)$ are homotopic.  
%where $\p C_r(c)=c\times \{r\}\subset c\times [0,r]/(c\times \{0\})$. 
As a result, $[f_{y, r}]=[f_{x, r}]\in H_2(X,X-K)$.
 Combining with  the observation from last paragraph we see 
 $[f_{y, r}]=[S(x)\cap \ol{B}(x,r)]\in H_2(X, X-K) $. 
 It follows that for any $y\in \ol{B}(x,r/4)$,
   $f_{y, r}$ represents a  nontrivial class in $H_2(X, X-\{y\})$ if and only if 
  $S(x)\cap \ol{B}(x,r)$ represents a nontrivial class in $H_2(X, X-\{y\})$.
Recall $S(x)\cap \ol{B}(x,r)$  is homeomorphic to the closed  unit  disk
 in the  Euclidean  plane.  Now it follows easily from excision that for any $y\in \ol{B}(x,r/4)$, 
    $S(x)\cap \ol{B}(x,r)$
%$ f_{y, r}$
 represents a nontrivial class in $H_2(X, X-\{y\})$ if and only if 
  $y\in \ol{B}(x,r/4)\cap (S(x)\cap \ol{B}(x,r))=S(x)\cap \ol{B}(x,r/4)$.
 %Therefore $supp(c)$ is a surface and for each
%$x\in supp(c)$ the circle $Link(supp(c),x)=Link(D,x)$ has length $2\pi$.  

\end{proof}

\begin{Cor}\label{s4}
   $supp(c)$ is locally convex in $X$.
\end{Cor} 

\begin{proof}
It follows from    the fact that 
for each
$x\in supp(c)$ the circle  $Link(supp(c), x)$  has length $2\pi$.

\end{proof}

An    argument  similar to the proof of    Lemma \ref{s3}  
 shows the complement of $supp(c)$ in $X-c$ is open:

\b{Le}\label{closed}
{$supp(c)$ is a closed subset of $X-c$.}
\end{Le}

Recall  $x_0$ lies in  some  open 2-cell  $A$.  
   Clearly   $A-x_0\xi\cup x_0\eta$  has two components. Let 
$A_0$ be the component of  $A-x_0\xi\cup x_0\eta$ 
that has    interior  angle less than $\pi$ at $x_0$. 

\begin{Le}\label{s5}
$A_0$  is contained in $supp(c)$.
\end{Le}

\begin{proof}
Let $y_0\in \p A_0-c$ be  the point  %on the boundary  of  $A_0$     
     with   %such that 
 $\angle_{x_0}(\xi,y_0)=\angle_{x_0}(y_0, \eta)$.   We claim if
  $x\not=x_0$ lies on  $x_0y_0$   and  $d(x_0, x)$   is  sufficiently small, then $x\in supp(c)$.

 Let   $x\in A_0$.  Then  the link $Link(X,x)$ is a circle with length $2\pi$.
Set $a_1=\log_x(x_0)$, $b_1=\log_x(\xi)$ and $c_1=\log_x(\eta)$.
Since $x\notin  x_0\xi\cup x_0\eta$, we have $d_x(a_1, b_1)<\pi$
 and  
$d_x(a_1, c_1)< \pi$.    On the other hand, 
$d_x(b_1, c_1)\leq d_T(\xi,\eta) < \pi$. 
When $x\ra x_0$ along  the segment $y_0x_0$, $x\xi\ra x_0\xi$ and $x\eta\ra x_0\eta$.
  It follows that for $x\in x_0y_0$ ($x\not=x_0$) with $d(x_0, x)$ sufficiently small, we have 
   $a_1\notin b_1c_1$, $b_1\notin a_1c_1$ and $c_1\notin a_1b_1$. 
    Thus   
 $Link(X,x)=a_1b_1\cup b_1c_1\cup c_1a_1$   for such  $x$. 
      The proof of Lemma \ref{s1} shows 
the path $\log_x(c)$ is homotopic to the path 
$ a_1b_1* b_1c_1* c_1a_1=Link(X,x)$.  Therefore    if
  $x\not=x_0$ lies on  $x_0y_0$   and  $d(x_0, x)$   is  sufficiently small, then
$\log_x(c)$  represents a nonzero class in $H_1(Link(X,x))$ and 
$x\in supp(c)$.

      So  we have  $A_0\cap supp(c)\not=\phi$.   Since   $A_0$ is disjoint from $X^{(1)}$,  
$supp(c)$ is a 2-dimensional manifold and 
   is closed in $X-c$,   we see $A_0\subset supp(c)$.

\end{proof}

Let $G$ be a finite metric graph and $L\subset G$  a homotopically nontrivial loop 
 in $G$. Then $L$ is homotopic to   a closed geodesic $l$ in $G$.  Since the universal cover of $G$ is 
  a simplicial tree and $L$ represents a hyperbolic isometry, we see that  $l\subset L$. 
   It follows that for any $x\in supp(c)$, the circle $c_x$ is contained in  $\log_x(c)$.

Fix a positive $r_0$ with $r_0 < r(x_0)$.  Then  $\partial \ol B(x_0,r_0)\cap A_0$
  is  an open arc contained in $A_0$. Fix two points 
$y, z\in \partial \ol B(x_0,r_0)\cap A_0$     with 
$\angle_{x_0}(y,z)\ge \angle_{x_0}(\xi, \eta)-\epsilon/4,$
and let $\sigma$ be the closed subarc of  $\partial \ol B(x_0,r_0)\cap A_0$ that joins 
$y$ and $z$.

\begin{Le}\label{s6}
For each  $p\in \sigma$,  there   exists  some  $\xi'\in \xi\eta$   with 
$p\in x_0\xi'$ and $x_0\xi'-\{x_0\}\subset supp(c)$. 
\end{Le}

\begin{proof} 
Fix $p\in \sigma\subset supp(c)$. We will try to extend the 
geodesic segment $x_0p$   inside $supp(c)$. 
Since $supp(c)$ is a surface and $Link(supp(c), p)$   is   a  circle with   length
$2\pi$, 
 the  geodesic segment $x_0p$ 
can be extended beyond $p$ in  $supp(c)$.  There are two 
cases to consider:  there is  either a finite  maximum
extension   %$x_0x$ with $p\in x_0x$  and $x\in X-supp(c)$, 
$$x_0x-\{x_0,x\}=(x_0p-\{x_0\})\cup (px-\{x\})\subset supp(c)\; \text{ where } \;x\in X-supp(c), $$
 %where $x\notin supp(c)$, 
or  an infinite maximal extension  % $x_0\xi'$ with $p\in x_0\xi'$   and  $\xi'\in \partial X$.    
$$x_0p-\{x_0\}\subset x_0\xi'-\{x_0\}\subset supp(c)\;\text{ where }\; \xi'\in \partial_\i X.$$

If the first case occurs and $x_0x-\{x_0,x\}$ ($x\notin supp(c)$) is  a finite maximal extension, 
 then Lemma \ref{closed} implies $x\in x_0\xi\cup x_0\eta$.  This is a    contradiction since 
 for any $x'$   in    $x_0\xi\cup x_0\eta$  
the geodesic segment $x_0x'\subset x_0\xi\cup x_0\eta$ does not 
pass through $p$.

Therefore the second case occurs and we have 
$x_0p-\{x_0\}\subset x_0\xi'-\{x_0\}\subset supp(c)$ 
  for some  $\xi'\in \partial_\i X$. 
 Pick a sequence of points $x_i$ on $ x_0 \xi'$ with $d(x_0, x_i)\rightarrow \infty$. 
  Since $c_{x_i}\subset \log_{x_i}(c)$  and $\log_{x_i}(\xi')\in Link(supp(c), x_i)=c_{x_i}$, 
    there is a point $\xi_i\in c$ with $\log_{x_i}(\xi')=\log_{x_i}(\xi_i)$. 
  $p\notin x_0\xi\cup x_0\eta$ implies $\xi_i\notin  x_0\xi\cup x_0\eta$. Therefore 
$\xi_i\in \xi\eta$. By the choice of $\xi_i$ the 
 sequence $\{x_0\xi_i\}_{i=1}^\i$   converges to the ray $x_0\xi'$. 
Since the sequence $\{\xi_i\}$ lies on the closed interval $\xi\eta$ we see
$\xi'\in \xi\eta$.  %The Lemma now follows. 

\end{proof}

Since for any $x\in supp(c)$ the circle $Link(supp(c), x)$ has length $2\pi$,     
   geodesics in $supp(c)$ do not branch.  It follows that 
for each   $p\in \sigma$,  there is a unique $\xi'\in \xi\eta$ 
  with  the property stated in Lemma  \ref{s6}. 
  Thus  we     can define a map $g:\sigma\rightarrow \xi\eta$ by $g(p)=\xi'$ where 
$\xi'$   is the unique point on $ \xi\eta$ with 
$x_0p-\{x_0\}\subset x_0\xi'-\{x_0\}\subset supp(c)$.

\begin{Le}\label{s7}
The map $g$ is continuous. 
\end{Le}

\begin{proof}
  The lemma follows easily from the facts that $\xi\eta$ is compact and   that 
%for any $x\in supp(c)$ the circle $Link(supp(c), x)$ has length $2\pi$ and thus 
geodesics in $supp(c)$ do not branch. 

\end{proof}

\begin{Le}\label{s8}
The map $\log_{x_0}$ is injective on 
the segment $\gamma(\epsilon/2)\gamma(h-\epsilon/2)$, and the surface     
$F_0=\rho_{x_0}(C( \gamma(\epsilon/2)\gamma(h-\epsilon/2)) )$  is convex in $X$.   
\end{Le}

% Let  $\xi_1=g(y)$, $\xi_2=g(z)$.
\begin{proof}
Let  $\xi_1=g(y)$, $\xi_2=g(z)$.   Notice $d_T(\xi_1, \xi_2)\ge  d_T(\xi,\eta)-\epsilon/2$
  holds by the     choice of $y$, $z$ and $x_0$.  
 %imply    $d_T(\xi_1, \xi_2)> d_T(\xi,\eta)-\epsilon/2$.  
We observe that it suffices to prove 
  the map $\log_{x_0}$ is injective on 
the segment $\xi_1\xi_2$ and the surface     
$\rho_{x_0}(C( \xi_1\xi_2) )$  is convex in $X$.

Since the map $g:\sigma\rightarrow \xi\eta$ is continuous and 
$\sigma$ is connected, the segment $\xi_1\xi_2$ is contained in 
$g(\sigma)$.  Therefore for each $\xi'\in \xi_1\xi_2$,    $x_0\xi'-\{x_0\}\subset supp(c)$  
%the open geodesic ray   $(x_0,\xi')$ lies in $supp(c)$ 
and $x_0\xi'\cap \sigma\not=\phi$.
Now the   map $\log_{x_0}$ is injective on 
%the segment 
$\xi_1\xi_2$
since  geodesics in $supp(c)$ do not branch.  The fact that
for each $x\in supp(c)$ the link $Link(supp(c),x)$ is a circle of length $2\pi$   
implies the surface     
$\rho_{x_0}(C( \xi_1\xi_2) )$  is   locally convex in $X$.  Since $\rho_{x_0}(C( \xi_1\xi_2) )$ 
  is also closed in $X$, it is convex in $X$.

\end{proof}

Lemma \ref{s8} completes the proofs of Theorems \ref{t16} and \ref{t13}.

\section{Applications}\label{applic}

Recall all  $CAT(0)$ $2$-complexes in this paper admit cocompact cellular isometric actions.  
In this section we give several applications of Theorem \ref{t13}. 
   But  first   we record  several   results concerning rays and flat sectors in 
a   $CAT(0)$  $2$-complex.  These 
     results  will be  used in  later    applications.

\subsection{Rays and Flat Sectors in    $CAT(0)$    $2$-complexes} \label{rays}

 %In this section we record  two  results concerning rays and flat sectors in 
%a CAT(0)  2-complex. The results will be used in the proof 
%of Theorem \ref{t15}.   
Recall flat sectors are defined in Section \ref{ideal}.  To simplify notation,  for any geodesic ray
 $\alpha:[0,\i)\ra X$,     we   also  use  $\alpha$   to   denote    its    image. 

%By a   \emph{sector}   we mean a closed 
%convex subset of $\mathbb{R}^2$ that is bounded  by two geodesic rays emanating from the origin.  A   sector is equipped with the
 %induced metric.     An isometric embedding from a sector  into a   CAT(0) space   $X$  is called a \emph{flat sector}  in $X$. 
%We notice   if $S$ is a flat sector in a  CAT(0) space   $X$, then 
%$\partial_T S$ is a closed interval  and isometrically embeds  into  the Tits boundary  $\partial_T X$.  

\begin{Prop}  \label{l1}
{Let   $X$ be a  $CAT(0)$    2-complex,     $c: [0,\i)\ra X$ a geodesic ray
    %that   admits a  cellular, isometric and cocompact  group action. 
    and       $S\subset X$  a flat sector. % and $c\subset  X$ a geodesic ray.
If $S\cap c=\phi$ and $c$ represents an interior point    of  $\partial_T S$,  then there is 
    a   geodesic  ray   $c':[0,\i)\ra  X$   asymptotic to $c$   
   with  $c'\subset S\cap X^{(1)}$,  
%$c'\subset S\cap X^{(1)}$ asymptotic to $c$,  
where    $X^{(1)}$  is the 1-skeleton of $X$. }
\end{Prop}
\begin{proof}  Since $X$ admits a cocompact    cellular isometric  action,  
%cocompact cellular and isometric  group action, 
 there is   some   $a>0$ with the following property:
  for any two 1-cells $vv_1$, $vv_2$ sharing an endpoint $v$,  if $\angle_v(v_1, v_2)<\pi$, then 
$d(v,v_1)+d(v,v_2)-d(v_1,v_2)>a$.   The flat sector $S$ is a closed convex subset of  $X$.
 Let  $\pi: X\rightarrow S$   be  the    orthogonal  projection.  Recall $\pi$ is 1-Lipschitz.
  %distance nonincreasing.
%The flat sector $S$ is a closed convex subset of  $X$, and  
 % the    orthogonal  projection   
%$\pi: X\rightarrow S$   is well-defined and 1-Lipschitz.  
 Let $p=c(0)$ and $q=\pi(p)$.  
Also let $p=p_1$, $p_2, \cdots, p_k=q$ be  a   finite sequence of points on the segment $pq$ with $d(p_i, p_{i+1})<{a}/{2}$.   
     Let  $c_i$    be  the  geodesic ray  with $c_i(0)=p_i$ and $c_i(\i)=c(\i)$.
%starting from $p_i$ asymptotic to $c$.  
 Notice if $c_i\cap S\not=\phi$ then $c_i([t_0, \infty))\subset S$
for some $t_0\ge 0$.    Since $c_1\cap S=c\cap S=\phi$ and $c_k\subset S$, there is 
some $i$ so that $c_i\cap S=\phi$ and   $c_{i+1}([t_0, \infty))\subset S$
for some $t_0\ge 0$.    Therefore it suffices to prove the lemma when $d(p,S)=d(p,q)<{a}/{2}$.

Now we assume  $d(p,S)=d(p,q)<{a}/{2}$.      As   $c$ represents a    point   in    $\partial_T S$,  
  the convex function $g(t):=d(c(t), S)$ is bounded from above, and is therefore non-increasing.  Thus 
$d(c(t), S)<{a}/{2}$  for all $t\ge 0$.   
The assumption  $S\cap c=\phi$    
implies 
for any $t$,  $\pi(c(t))$   lies either  on the   manifold boundary of   $S$ or 
  in $S\cap X^{(1)}$.    Since  $c$ represents an interior point    of  $\partial_T S$,   for $t$ sufficiently large 
  $\pi(c(t))$  has to lie in 
the interior of $S$ and thus in $S\cap X^{(1)}$.  We may assume $\pi(c(t))\in S\cap X^{(1)}$ for 
all $t$ by considering a sub-ray of $c$.   
 Then $\pi(c)\subset  S\cap X^{(1)}$ is  a continuous path  going to infinity.  
%contained in   $S\cap X^{(1)}$.    
Let $t_1< t_2<t_3$ so that $v_1:=\pi(c(t_1))$, $v_2:=\pi(c(t_2))$,  $v_3:=\pi(c(t_3))$ 
are 0-cells,    and   $v_1v_2$ and  $v_2v_3$ are  1-cells.   %     and  $v_2$, $v_3$ are  also joined by a 1-cell.   
We claim $\angle_{v_2}(v_1, v_3)=\pi$.   The claim clearly implies   there is 
a geodesic ray    $c':[0,\i)\ra  X$  asymptotic to $c$    with
$c'\subset \pi(c)  \subset  S\cap X^{(1)}   $.
%asymptotic to $c$.  
 Suppose  $\angle_{v_2}(v_1, v_3)< \pi$.   The choice of  $a$ implies $a+d(v_1, v_3)<d(v_1, v_2)+d(v_2, v_3)$.
 Triangle inequality implies 

\begin{tabular}{ll}
     & $d(c(t_1), c(t_3))$\\
$\le$  & $    d(c(t_1), \pi(c(t_1)))+d(\pi(c(t_1)),  \pi(c(t_3)))+ d( \pi(c(t_3)), c(t_3))$\\
$\le$   &   ${a}/{2} +d(\pi(c(t_1)),  \pi(c(t_3)))+  {a}/{2}$\\
%$ =$   &   $a+d(\pi(c(t_1)),  \pi(c(t_3)))$\\
  $=$   &   $a+d(v_1, v_3)$\\
$<$   &   $d(v_1, v_2)+d(v_2, v_3),$
\end{tabular}\newline
or  $d(c(t_1), c(t_3))< d(v_1, v_2)+d(v_2, v_3)$.  
  On the other hand,  since $\pi$ is 1-Lipschitz, the   length  of $ \pi\circ c_{|[t_1, t_3]}$ 
 is less than or equal to the   length   of $c_{|[t_1, t_3]}$.   It follows   

\begin{tabular}{ll}
     & $d(v_1, v_2)+d(v_2, v_3)$\\
$=$    &    $d(\pi(c(t_1)),  \pi(c(t_2)))+ d( \pi(c(t_2)), \pi(c(t_3)))$\\
$\le$    &     $length(\pi\circ c_{|[t_1, t_2]})+ length(\pi\circ c_{|[t_2, t_3]})$\\
$=$    &   $length(\pi\circ c_{|[t_1, t_3]})$\\
$\le$   &    $length(c_{|[t_1, t_3]})$  \\
$=$   &   $d(c(t_1), c(t_3))$, 
\end{tabular}  \newline
or  $d(v_1, v_2)+d(v_2, v_3)\le  d(c(t_1), c(t_3))$.  The contradiction proves the claim.   

  \end{proof}

\begin{Cor}\label{c2}
{Let     $X$    be   a      $CAT(0)$    2-complex
 %that    admits a  cellular, isometric and cocompact  group action,   
 and       $S_1,  S_2\subset  X$ 
   two flat sectors. 
  If $\partial_T S_1 \cap\partial_T S_2$ is a nontrivial interval, then 
$S_1\cap S_2\not=\phi$.  } 
\end{Cor}
\begin{proof}   Suppose  $S_1\cap S_2=\phi$.   For each $\xi\in interior(\partial_T S_1 \cap\partial_T S_2)$, 
      let   $c_{\xi}$   be  a   ray  in $S_2$ asymptotic to $\xi$.   Then  $c_{\xi}\cap S_1=\phi$.     
   Proposition    \ref{l1} 
     implies there is a ray
  $c'_{\xi}\subset S_1\cap X^{(1)}$ representing $\xi$.   There are uncountably many points in 
 $interior(\partial_T S_1 \cap\partial_T S_2)$, so there are  uncountably many rays in 
 $S_1\cap X^{(1)}$   pointing to    uncountably many  directions.  On the other hand, there are only
 countably many 1-cells in  $S_1\cap X^{(1)}$ and they can  only  give rise     to  countably many directions. 
 The contradiction proves the corollary.   

 \end{proof}

For  a  flat sector $S\subset X$ and $x\in S$, let $S(x)\subset S$ be 
  the   flat sector with cone point $x$
  and  $\p_T S(x)=\p_T S$.

\begin{Prop}\label{secto}
{Let  $X$ be a  $CAT(0)$   $2$-complex  and $\sigma:[0, h]\ra \p_T X$ a geodesic segment with length
  $h\le \pi$.  Suppose there are numbers $t_i$ $(1\le i \le 4)$  such that  $0<t_1<t_2<t_3<t_4< h$,
   and two flat sectors $S_1$ and $S_2$ with $\p_T S_1=\sigma([t_1,t_3])$ and 
$\p_T S_2=\sigma([t_2,t_4])$. Let $c_1: [0,\i)\ra S_1$    and  $c_2: [0,\i)\ra S_2$ 
   be  rays   with $c_1(+\i)=\sigma(t_1')$,
  $t_1<t_1'<t_2$    and   $c_2(+\i)=\sigma(t_2')$,   $t_3<t_2'<t_4$  respectively.  Then   
  there    is some $u_0\ge 0$ such that $S_1(c_1(u))\cap S_2(c_2(u'))$ is a flat sector
   for all $u, u'\ge u_0$.}
\end{Prop}

 The   proof of Proposition \ref{secto}
 is divided into a   few lemmas.

\begin{Le}\label{outg}
{There is   some $a\ge 0$ such that
 $c_1([a, \i))\cap S_2=\phi$ and $c_2([a, \i))\cap S_1=\phi$.}
\end{Le}

\begin{proof}  Assume    $c_1\cap S_2\not=\phi$. 
  Since  $c_1\cap S_2$ is a closed convex subset  of $c_1$,  it is 
either a closed segment
  of   $c_1$  or a subray of $c_1$.   But $c_1(+\i)=\sigma(t_1')\notin \p_T S_2$  implies 
  $c_1\cap S_2$   can not  be a  ray.   Therefore    
 there is some $a\ge 0$ such that  $c_1([a, \i))\cap S_2=\phi$.  The proof 
 of the  second equality   is similar.

\end{proof}

\begin{Le}\label{theri}
{With notation   as   in  Proposition  \ref{secto}  and  Lemma \ref{outg}.    
   Then there are $u_1, u_2\ge a$ such that
  $S_1(c_1(u_1))\cap S_2(c_2(u_2))$ is a flat sector  
 with cone point   $c_1(u_1)\sigma(t_3)\cap c_2(u_2)\sigma(t_2)$.}

%and the two rays 
  %bounding $S_0$ lie in  $c_1(u_1)\sigma(t_3)$ and $c_2(u_2)\sigma(t_2)$  respectively. }
\end{Le}

\begin{proof}
Since $\p_T S_1\cap \p_T S_2=\sigma([t_2, t_3])$ is a nontrivial interval,
    Corollary \ref{c2} implies $S_1\cap S_2\not=\phi$.  $S_1\cap S_2$ is a closed convex subset 
  of $S_i$ ($i=1,2$). For any point $p\in S_1\cap S_2$, let  $S_p$ be the   flat  sector 
 with cone point $p$ and   $\p_T S_p=\sigma([t_2, t_3])$.  Clearly
  we have     $S_p\subset S_1\cap S_2$.  Fix a point $p\in S_1\cap S_2$ and  
   consider the  subsets $S_p$ and $c_2$ of $S_2$.
  It    is   clear that   for large enough $u$ we have   $c_2(u)\sigma(t_2)\cap S_p\not=\phi$.
  It follows that  $c_2(u)\sigma(t_2)\cap S_1\not=\phi$  
for large enough $u$.   Fix a $u_2\ge a$ with  $c_2(u_2)\sigma(t_2)\cap S_1\not=\phi$.
     Note   $c_2(u_2)\sigma(t_2)\cap S_1$ is  a  subray of $c_2(u_2)\sigma(t_2)$.

  Consider the   rays $c_2(u_2)\sigma(t_2)\cap S_1$  and $c_1$   in  $S_1$. 
  Since $t_1<t_1'<t_2<t_3$    %, $\p_T S_1=\sigma([t_1, t_3])$, 
      and  $c_1(+\i)=\sigma(t_1')$,
  we see for large enough $u$, $c_1(u)\sigma(t_3)\cap c_2(u_2)\sigma(t_2)\not=\phi$.
  %Clearly $c_1(u)\sigma(t_3)\cap c_2(u_2)\sigma(t_2)$ consists of a point. 
   Fix a   $u_1\ge a$ such that $c_1(u_1)\sigma(t_3)\cap c_2(u_2)\sigma(t_2)\not=\phi$.
  Clearly $c_1(u_1)\sigma(t_3)\cap c_2(u_2)\sigma(t_2)$  consists of a    single  point.
  Let $x$ be this point.     
   We notice  $x\in S_1\cap S_2$  and   that  $S_x\subset S:=S_1(c_1(u_1))\cap S_2(c_2(u_2))$.  We claim 
$S_x=S$.

        Suppose $S_x\not=S$.  Pick $y\in S-S_x$.     Notice 
  $\{y\}$ and $S_x$ are contained in the flat sector    $S_1(c_1(u_1))$.
   Since $x$ lies on the boundary ray $c_1(u_1)\sigma(t_3)$  of  $S_1(c_1(u_1))$,
      $0< \angle_x(y, \sigma(t_2))<   \angle_x(y, \sigma(t_3))\le \pi$.
  Similarly  by   viewing  $\{y\}$ and $S_x$   as subsets of $S_2(c_2(u_2))$  we 
   see  $0<\angle_x(y, \sigma(t_3))< \angle_x(y, \sigma(t_2))\le \pi$.
   A contradiction.

\end{proof}

\begin{Le}\label{somes}
{With notation   as   in  Proposition  \ref{secto}  and  Lemma \ref{outg}. If  
   there are $u_1, u_2\ge a$ such that
  $S_1(c_1(u_1))\cap S_2(c_2(u_2))$ is a flat sector  with  
$c_1(u_1)\sigma(t_3)\cap c_2(u_2)\sigma(t_2)$   as    cone point, 
then for all 
$u\ge u_1$, $u'\ge u_2$,   
$S_1(c_1(u))\cap S_2(c_2(u'))$ is a flat sector.}
%with cone point $c_1(u)\sigma(t_3)\cap c_2(u')\sigma(t_2)$.}
\end{Le}

\begin{proof}
Let $S=S_1(c_1(u_1))\cap S_2(c_2(u_2))$.  Then   $\p_T S=\sigma([t_2, t_3])$. 
   Notice both $S$ and 
 $S_2(c_2(u'))$ ($u'\ge u_2$)  
are subsectors of the flat sector $S_2(c_2(u_2))$.  It   is   clear that
  $S\cap S_2(c_2(u'))$  is  a flat sector  with cone point  $c_1(u_1)\sigma(t_3)\cap c_2(u')\sigma(t_2)$.
   It follows   from   
  $S_1(c_1(u_1))\cap  S_2(c_2(u'))\subset S_1(c_1(u_1))\cap S_2(c_2(u_2))\cap S_2(c_2(u'))
=S\cap S_2(c_2(u'))$    that   
$S_1(c_1(u_1))\cap  S_2(c_2(u'))=S\cap S_2(c_2(u'))$
is  a   flat sector  with 
$c_1(u_1)\sigma(t_3)\cap c_2(u')\sigma(t_2)$       as  cone point.    
    Now a  similar   argument, but 
fixing $u'$ and increasing  $u_1$, shows that  for all 
$u\ge u_1$,
  $S_1(c_1(u))\cap S_2(c_2(u'))$ is a flat sector.

\end{proof}

The   proof of Proposition \ref{secto}
 is now complete.

\subsection{Circles  in Tits Boundary}

  %  $\equiv$  $\not\equiv$ 

We recall  an   \e{$n$-flat}  in a   $CAT(0)$   space   $X$  is the image of an isometric embedding from the $n$-dimensional 
Euclidean space into  $X$.

\begin{Th}\label{t14}\emph{(V. Schroeder \cite{BGS},    B. Leeb \cite{L})}
{Let $X$ be a locally compact  $CAT(0)$ space. Suppose $\mathbb{S}^{n-1}\subset \partial_TX$ is a 
unit $(n-1) $ sphere in the Tits boundary that   does not bound a unit hemisphere. Then 
   there is an   $n$-flat $F\subset X$ such that 
$\partial_T  F=\mathbb{S}^{n-1}$. }
\end{Th}

 Let  $X$ be a  $CAT(0)$   $2$-complex.  %    with   a cellular, isometric and cocompact group  action. 
    Theorem \ref{t14}   and Theorem \ref{t11}    imply  that   any unit circle in   $\p_T X$  
is the ideal boundary of  a 2-flat.
Theorem \ref{t13}  enables us to generalize  this   result    to topological circles in the Tits boundary. 
  Recall  Corollary  \ref{r1} implies topological
circles in $\p_T X$ 
 are  simple closed geodesics.

%\begin{Le}\label{flatse}
%{Let  $X$ be a  $CAT(0)$   $2$-complex  and $\sigma:[0, h]\ra \p_T X$ a geodesic segment with length
 % $h\le \pi$.  Suppose there are numbers $t_i$ ($1\le i \le 4$) with $0<t_1<t_2<t_3<t_4< h$,
  % and two flat sectors $S_1$ and $S_2$ with $\p_T S_1=\sigma([t_1,t_3])$ and 
%$\p_T S_2=\sigma([t_2,t_4])$.  Given any $p\in S_1\cap S_2$,   let $c_3\subset S_1$ ($c_2\subset S_2$)
 %  be the ray   that contains the ray $p\sigma(t_3)$ ($p\sigma(t_2)$) and   separates the sector
  %$S_1$ ($S_2$) into two components. %Denote by  $S$ the unique sector with cone point $p$ and
 %$\p_T S=\sigma([t_2, t_3])$.   
 %Let  $\tilde S_1$ ($\tilde S_2$)   be a   flat sector
 %with cone point  on $c_3$ ($c_2$)  and $\p_T {\tilde S_1}=\p_T S_1$ ($\p_T {\tilde S_2}=\p_T S_2$ ). 
  %Then    $S:=\tilde S_1\cap \tilde S_2$  is  a  flat sector with 
   %$\p_T S=\sigma([t_2, t_3])$.}

Let $A\subset X$  be  a subset of a $CAT(0)$ space.  A point $\xi\in \p_\i X$ is a \e{limit point}
of $A$ if there is a sequence $a_i\in A$ ($i\ge 1$) such that $\{a_i\}$ converges to
$\xi$ in the cone topology. 
%For any $CAT(0)$ space $X$ and any subset $A\subset X$, let $\ol{A}\subset \ol{X}$  be the closure
 %of $A$  in $\ol{X}$ with respect to the  cone topology.  
The \e{limit set}   $L(A)\subset \p_\i X$
  of $A$ is  the set of limit points of $A$. 
Recall a quasi-flat in     a  $CAT(0)$ $2$-complex    $X$ is the image of a quasi-isometric 
  embedding from the Euclidean plane $\mathbb{R}^2$ into    $X$.  
 
\begin{Prop}\label{c3}
{ Let $X$ be a  $CAT(0)$   $2$-complex %    with   a cellular, isometric and cocompact group  action
        and  
$C\subset \partial_TX$  a topological circle in the Tits boundary. Then there is a quasi-flat  $E$ of $X$ with the following properties:\newline
\e{(i)}    $E$ is homeomorphic to the plane with the closed unit disk removed;\newline
\e{(ii)}    $E$ is flat, i.e.,    each point of $E$ has a neighborhood in $E$ which is isometric to an open subset of $\mathbb{R}^2$;\newline
\e{(iii)}    $L(E)=C$.}
\end{Prop}
\begin{proof}   By Corollary \ref{r1}  we may assume the circle $C$ is a simple closed geodesic.    
Choose points $a_1,a_2,\cdots, a_n$ on $C$ in   cyclic  order such that they divide $C$ into intervals of 
equal length $l< \pi/4$.   For each $i$ (mod $n$)  let 
 $m_i$  be the  midpoint of $a_ia_{i+1}$.   
 Theorem \ref{t13}  implies that  there are flat sectors $S_i$ ($1\le i\le n$) such that
 $a_ia_{i+1}$ is contained in the interior of $\p_T S_i$ and $m_j\notin \p_T S_i$  for $j\not=i$.
  Let $c_i:[0,\i)\ra S_i$ be a ray in $S_i$ with $c_i(+\i)=m_i$. 
   Now  Proposition \ref{secto}  implies that there is some $u_0\ge 0$ such that 
for any $u_i\ge u_0$ ($1\le i\le n$), the intersection 
 $S_i(c_i(u_i))\cap S_{i+1}(c_{i+1}(u_{i+1}))$ ($i$ mod $n$)
  is a flat sector.   Since  $\p_T S_i\cap \p_T S_j=\phi$ 
for any $1\le i, j\le n$  with  $i-j\not\equiv -1,0,1 $ mod $n$,  there is some  $u'_0\ge 0$
   such that for any $u_i\ge u'_0$ ($1\le i\le n$), the intersection 
 $S_i(c_i(u_i))\cap S_{j}(c_{j}(u_{j}))=\phi$ 
for  any $1\le i, j\le n$  with   $i-j\not\equiv -1,0,1 $ mod $n$. Now choose
  $u_i\ge u_0, u'_0$   and 
   let $E_i$ be the interior of the flat sector   $S_i(c_i(u_i))$.
     Set   $E=\cup_i   E_i$.  Now   it  is   easy to see 
   that $E$ is a quasi-flat   with   the desired properties.

 \end{proof}

We  call  a geodesic $c$ of the form $c: (0,\i)\ra X$ in a $CAT(0)$ space  $X$
an \e{open geodesic ray}.

\begin{Cor}\label{r2}
{Let $C\subset \p_T X$ and $E\subset X$ be  as in Proposition \ref{c3}.   Then $E$ admits a foliation by   open geodesic rays 
   with  the following  properties:\newline
\e{(i)} each ray in the foliation is asymptotic to a point in $C$,  and   for each $\xi\in C$ 
  there is at least one ray in the foliation asymptotic to $\xi$;\newline
\e{(ii)} there is a constant $a>0$ such that the distance between    any  
    two asymptotic rays 
 in the foliation is at most $a$.}
\end{Cor}

\begin{proof}
  We use the  notation  in the proof of Proposition \ref{c3}. 
  Let $p_i$  denote  the cone point of the flat sector 
$S_{i-1}(c_{i-1}(u_{i-1}))\cap  S_i(c_i(u_i))$. 
  For each $i$ ($1\le i\le n$), let $S_i'\subset \ol{E}$ be the  flat sector with 
   cone point  $p_i$   and   
    $\p_T S_i'=m_{i-1}m_i$.   % and  having   the same cone point 
  %as that of $S_{i-1}(c_{i-1}(u_{i-1}))\cap  S_i(c_i(u_i))$.  
  As a flat sector $S_i'$ is  a   union of rays issuing from its cone point. 
    $S_i(c_i(u_i))-S_i'\cup S_{i+1}'$ is basically a flat strip:
 it is a convex subset of the flat sector   $S_i(c_i(u_i))$  whose boundary contains two asymptotic rays.
  We foliate  $S_i(c_i(u_i))-S_i'\cup S_{i+1}'$
by parallel rays. In this way 
 we get a foliation  of $E$ by open geodesic rays  with the desired properties.

\end{proof}

One    question about circles in the Tits boundary is whether the lengths of circles in the Tits boundary
form a discrete set.   % In the special case when  all the interior angles of 2-cells are rational,  
  %the following result    provides a positive answer. 

 \begin{Th}\label{circlelength}
{Let $X$ be a   $CAT(0)$     $2$-complex.
 %that   admits a  cellular, isometric and cocompact   action. 
  If the interior angles of all the 2-cells of $X$ are 
rational multiples of $\pi$, then there is a positive integer $m$ so that the length of each topological circle in the Tits boundary is an integral  multiple of $\pi/m$. }
\end{Th}

The following lemma is not hard to prove. 
%This theorem follows from Proposition \ref{c3}  and the following easy lemma.

\begin{Le} \label{polygon}
{Let $p_1$, $p_2, \cdots,  p_n\in {\mathbb R}^2$   be $n $ points in the  Euclidean plane, $a$, $b$  be  two rays emanating from 
$p_1$, $p_n$ respectively,  and $e_i$ $  (i=1, 2 , \cdots, n-1)$  the segment
   connecting $p_i$ and $p_{i+1}$.  Suppose the union 
$c:=a\cup e_1 \cup  \cdots \cup e_{n-1}\cup b$ is a simple path and $A$ is one of the two 
  components of ${\mathbb R}^2-c$.
 %regions with boundary  $c$. 
   % Let $\partial_T A\subset \p_T {\mathbb R}^2$  be the segment that is the \lq\lq ideal boundary" of $A$.
Denote the interior angle of $A$ at $p_i$ by $A_i$. Then the length of the   
closed interval  $L(A)\subset \p_T {\mathbb R}^2$  is 
$\Sigma_{i=1}^{n}A_i- (n-1)\pi$.}
\end{Le}            
 \qed

\vspace{5mm}

\noindent
{\bf{Proof of Theorem \ref{circlelength}:}}

Since $X$ admits a cocompact cellular isometric action, 
 the  assumption implies there is a positive integer $m$ so that  the interior angles  of 2-cells are integral  multiples of 
$\pi/m$. 

Let $C$ be a topological circle in     $\p_T X$  and $E$   a  quasi-flat provided by  Proposition  \ref{c3}
with  $L(E)=C$.  Corollary \ref{r2} implies $E$ admits  a    foliation by open geodesic
rays.   Choose points $a_1,a_2,\cdots, a_n$ on $C$ in   cyclic  order such that they divide $C$ into intervals 
  with  length $< \pi/2$.   For each $i$ ($1\le i\le n$), let $c_i$ be an open geodesic ray 
  in the foliation of $E$ that represents $a_i$.   Then $c_i\cap X^{(1)}\not=\phi$.  
  Pick a point $x_i\in c_i\cap X^{(1)}$  as follows. If $c_i\cap X^{(1)}$  contains vertices then let
 $x_i$ be a fixed vertex in $c_i\cap X^{(1)}$, otherwise if $c_i\cap X^{(1)}$  contains no  vertex
  let  $x_i$ be any  fixed point in $c_i\cap X^{(1)}$.  

Notice   $E-\cup c_i$  has 
  %The $n$   open geodesic rays  $c_i$ divides $E$ into 
$n$ components, each of which is  contractible.
  Let $E_i$ be the   component of  $E-\cup c_i$    that   contains  
 $c_i$ and $c_{i+1}$ ($i$ mod $n$) on the boundary. 
   Let $\alpha_i\subset E\cap X^{(1)}$ be a    simple path  from  $x_i$ to  $x_{i+1}$ 
  with   $interior(\alpha_i)  \subset E_i$.  Then $\alpha_i$ divides $E_i$ into two components.
  Let $D_i$ be the  unbounded component of $E_i-\alpha_i$,  $A_{i,j}$ ($j=2,\cdots, k_i-1$)  the 
 interior angles of $D_i$ at vertices in the interior of the   path  $\alpha_i$,  and 
 $A_{i,1}$  and  $A_{i, k_i}$  the 
 interior angles of $D_i$ at $x_i$    and    $x_{i+1}$  respectively. 
   %at the vertices on the path $\alpha_i$ from   $x_i$ to  $x_{i+1}$.  
  Notice   $A_{i,j}$ ($1<j<k_i$)  and $A_{i+1,1}+A_{i, k_i}$ ($i$ mod $n$)
     are    integral  multiples of 
$\pi/m$.   
Now the theorem follows by applying Lemma \ref{polygon}  to each $D_i$ and adding up all
  the equalities.

% Recall $E$ is homeomorphic to the plane with the unit disk removed.  Choose a simple closed path $c\subset E\cap X^{(1)}$  so that the inclusion $c\subset E$ is a homotopy equivalence,  and let $E'$  be  the unbounded component 
%of $E-c$.   We decompose $E'$ into a union of subsets that have the same form as the $A$ in Lemma \ref{polygon}
 %   by extending 
%certain edges in $c$ to rays. Now the theorem follows by applying Lemma \ref{polygon}  to each of these subsets and observing that 
%all the angles are integer multiples    of $\pi/m$.  

\qed

\subsection{   Branch Points  in Tits Boundary}

   Let $X$ be a  $CAT(0)$    $2$-complex.
    %that    admits   a cellular, isometric and cocompact group  action.
       %Corollary \ref{c1} 
  Since  small metric balls in    $\partial_T X$  are $R$-trees  (Corollary \ref{c1})
    we can talk about branch points in $\p_T X$ (see Definition \ref{branchpt}). 
   The following proposition says  a branch point  in $\p_T X$  is represented 
  by a geodesic ray where flat sectors   branch off.

\begin{Prop}\label{branchpo}
{Let $X$ be a   $CAT(0)$      $2$-complex
   %that    admits a  cellular, isometric and cocompact  group action, 
and $\xi\in \p_T X$ a 
  branch point.  Then there are   two flat sectors $S_1$ and $S_2$   and a 
 geodesic ray $c:[0,\i)\ra X$  such that \newline
\e{(i)} $c(0)$ is the common cone point of $S_1$ and $S_2$;\newline
\e{(ii)} $c(+\i)=\xi$   lies in the  interior   of the segment  $\p_T S_i$  $(i=1,2)$;\newline
\e{(iii)} $S_1\cup S_2-c$ has three  components and the closure of 
 each component  is a flat sector; in particular $c\subset X^{(1)}$.}
 
\end{Prop}

\begin{proof}
  Since  $\xi\in \p_T X$ is a branch point,  there are points 
$\xi_1,  \xi_2, \xi_3\in \p_T X$      such  that    $0<d_T(\xi, \xi_i)<\pi/4$ and 
$\xi_1\xi_2\cap \xi_1\xi_3=\xi_1\xi$.   By   Theorem \ref{t13}   we may assume 
  there are flat sectors $S'_1$ and $S'_2$ such that 
 $\p_T S'_1=\xi_1\xi_2$,  $\p_T S'_2=\xi_1\xi_3$.  Corollary \ref{c2} implies 
    $S'_1\cap S'_2\not=\phi$.    
Thus $S'_1\cap S'_2$   is a nonempty closed convex subset of the flat sector $S_1'$. 
  It follows that $\p_T(S'_1\cap S'_2)=\xi_1\xi$.  Notice
  for any $x$ in the interior of $S'_1\cap S'_2$, the intersection
  $x\xi_2\cap  (S'_1\cap S'_2)$ is a  segment.  
  Since  $x\xi_2$ lies in the interior of $S_1'$,  
    the   intersection of $x\xi_2$   with 
  the boundary of $S'_1\cap S'_2$   lies   either in $X^{(1)}$   or  on $o_2\xi_3$,  where 
  $o_2$ denotes the cone point of $S_2'$.
 Since  $\xi_3\notin  \p_T(S'_1\cap S'_2)$,   the intersection of $x\xi_2$   with 
  the boundary of $S'_1\cap S'_2$  can not lie on $o_2\xi_3$  when  $d(o_2, x)$ is 
  sufficiently large.    So there is an infinite path in the boundary of 
  $S'_1\cap S'_2$  that is contained in  $X^{(1)}$. 
The fact that $X$ admits a cocompact 
cellular isometric action  implies there is a geodesic ray $c\subset X^{(1)}$  such that 
  $c$ lies in the boundary of $S'_1\cap S'_2$.  It follows    $c(+\i)=\xi$.
%for any $x$ in the interior of $S'_1\cap S'_2$,
   %$x\xi_2\cap c\not=\phi$.

  Set $S_1=S'_1(c(0))$, $S_2=S'_2(c(0))$.  Also  let $S\subset S_1\cap S_2$ be the flat sector 
  with cone point $c(0)$ and $\p_T S=\xi_1\xi$.    The choice of $c$ implies that 
   $S_1$  and $S_2$ branch off at $c$  and $S=S_1\cap S_2$.

\end{proof}

%An $R$-tree may have lots of branch points ( see Definition \ref{branchpt}).     
We would like to know 
   if  there   is   a   positive 
constant $c=c(X)$ such that the distance between any two branch points in $\partial_T X$
 is at least $c$.    
 The positive answer to the question would imply that the   components of the 
    Tits boundary are   almost like   simplicial metric graphs. Recall that the Tits boundary of a 2-dimensional Euclidean building is a 1-dimensional spherical building.

 \begin{Th}\label{t15}
{Let $X$ be a   $CAT(0)$      $2$-complex.
%   that    admits a  cellular, isometric and cocompact  group action. 
If the interior angles of all the 2-cells of $X$ are 
rational multiples of $\pi$, then there is a positive integer $m$   
  such that    the distance between any two branch points in $\partial_T X$ is 
  either infinite or   an integral  multiple of 
      ${\pi}/{m}$.}
\end{Th}
 
\begin{proof}  Let $\xi$ and $\eta$ be two branch points in $\partial_T X$ with $d_T(\xi,\eta)< \i$.
  Let   $\sigma\subset \p_T X$   be  a minimal  geodesic from $\xi$ to $\eta$. 
   Since $\xi$  is  a   branch point, there are points $\xi_i$ ($i=1,2,3$) in $\p_T X$  with  $\xi_3\in \sigma$   and 
  $0<d_T(\xi_i, \xi)<\pi/4$  such that $\xi_3\xi_1\cap \xi_3\xi_2=\xi_3\xi$.  Similarly
   there are points $\eta_i\in \p_T X$ ($i=1,2,3$)  with  $\eta_3\in \sigma$   and 
  $0<d_T(\eta_i, \eta)<\pi/4$  such that $\eta_3\eta_1\cap \eta_3\eta_2=\eta_3\eta$. 
Set $\sigma_1=\xi_1\xi\cup \sigma\cup \eta\eta_1$ and 
$\sigma_2=\xi_2\xi\cup \sigma\cup \eta\eta_2$.   

By Theorem \ref{t13}  and Proposition \ref{secto} (see the proof of Proposition \ref{c3})
      there are flat sectors  $S_1, \cdots, S_k$ such that
 $S_i\cap S_{i+1}$  is a flat sector, $S_i\cap S_j=\phi$ for $|j-i|\ge 2$   and  $L(E_1)=\sigma_1$,  where $E_1=\cup_i {S_i}$.  We can foliate $E_1$  by open geodesic rays as in Corollary  \ref{r2}. 
  Similarly  there are flat sectors  $S'_1, \cdots, S'_k$ such that
 $S'_i\cap S'_{i+1}$  is a flat sector,  $S'_i\cap S'_j=\phi$ for $|j-i|\ge 2$   
    and  $L(E_2)=\sigma_2$,  where $E_2=\cup_i {S'_i}$.  
  The proof of Proposition  \ref{branchpo}
  shows  that there are geodesic rays  $c'_1, c'_2\subset E_1\cap X^{(1)}$
   such that $c'_1(+\i)=\xi$, $c'_2(+\i)=\eta$.

  Notice  the assumption implies there is a positive 
integer $m$ so that  the interior angles  of 2-cells are integral  multiples of 
$\pi/m$.   Choose points $\xi=a_1,a_2,\cdots, a_n=\eta$ on $\sigma$ in   linear   order such that 
they divide $\sigma$ into intervals 
  with  length $< \pi/2$.  Set  $c_1=c'_1$ and $c_n=c'_2$. 
 For each $i$ ($1< i< n$), let $c_i$ be an open geodesic ray 
  in the foliation of $E_1$ that represents $a_i$.   
     Now  Lemma   \ref{polygon}
  and the proof of Theorem  \ref{circlelength}  show that 
  $d_T(\xi, \eta)$ is   an integral  multiple of   $\pi/m$.

 \end{proof}

\subsection{$\pi$-Visibility}  \label{vis}

Let  $X$ be a   locally   compact $CAT(0)$ space. 
      A \emph{flat half-plane}  in $X$  is 
   the  image of an isometric embedding  
  $f:\{(x,y)\in \E^2: y\ge 0\}\rightarrow X$, and in this case we say the geodesic 
  $c:R\rightarrow X$,  
$c(t)=f(t,0)$   bounds the 
flat half-plane.   

Let    $\xi, \eta\in \p_T X$ with $d_T(\xi,\eta)=\pi$.  
If there is a geodesic  $c$ in $X$   with $\xi$    and  $\eta$   as   endpoints,   then $c$ bounds a flat 
 half-plane (see \cite{BH}).  
In general,    there is  no geodesic in $X$ with 
$\xi$    and  $\eta$   as   endpoints. 

 Recall  if  $X$ is  a    $CAT(0)$ 2-complex,   %   that   admits an isometric cocompact  group action,  
 then for any  $\xi\in \partial_T X  $ and any $r: 0 < r < {\pi}/{2}$,
      the  closed metric ball $\ol B(\xi, r)$ is an $R$-tree.  It is necessary to
 distinguish those points in  $\p_T X$ that are \lq\lq dead ends''.

\b{Def}\label{terminalp}
{Let $X$ be a $CAT(0)$ 2-complex  and $\xi\in \p_T X$.  We
     say $\xi$ is a \e{terminal point}
 if $\xi$ does not lie in the interior of any geodesic segment in $\p_T X$.}
\end{Def}

\b{Th}\label{pivisi}
{Let   $X$ be a  $CAT(0)$  2-complex.
    %that   admits a  cellular, isometric and cocompact  group action. 
%a  proper, cocompact   action by cellular isometries.
If $\xi,\eta\in \p_T X$ are not terminal points and $d_T(\xi,\eta)\ge \pi$, then there is a geodesic
 in $X$   with $\xi$    and  $\eta$   as   endpoints.}
\end{Th}

\b{proof}   By Proposition \ref{p2} 
    we may assume $d_T(\xi,\eta)=\pi$.  Let $\sigma:[0,\pi]\ra \p_T X$ be a 
 minimal geodesic from $\xi$ to $\eta$.  Since $\xi,\eta\in \p_T X$ are not terminal points 
 and small metric balls in   $\p_T X$  are $R$-trees, there is some $\epsilon$, $0<\epsilon<\pi/4$ such that 
  $\sigma$  extends to 
 a locally isometric map  $[-\epsilon, \pi+\epsilon]\ra \p_T X$, which  we  still 
   denote   by   $\sigma$.

Notice $\sigma_{|[-\epsilon, \pi/2+\epsilon]}$ is a minimal geodesic   in  $\p_T X$
with length less than 
 $\pi$.     By Theorem \ref{t13} there is  a flat sector $S_1$ in $X$ such 
that $\p_T {S_1}=\sigma([-\epsilon/2, \pi/2+\epsilon/2])$.    Similarly there is a 
 flat sector $S_2$ in $X$ with $\p_T {S_2}=\sigma([\pi/2-\epsilon/2, \pi+\epsilon/2])$. 
 Since $\p_T {S_1}\cap \p_T {S_2}=\sigma([\pi/2-\epsilon/2, \pi/2+\epsilon/2])$
is a     nontrivial interval,    Corollary \ref{c2}  implies 
%$X$ is a 2-complex and the intersection of the boundaries of $S_1$ and $S_2$ is a
% nontrivial interval,   we have 
 $S_1\cap S_2\not=\phi$.  Pick  $x\in S_1\cap S_2$ and let
 $S\subset S_1\cap S_2$ be the subsector with cone point $x$ and 
 $\p_T S=\sigma([\pi/2-\epsilon/2, \pi/2+\epsilon/2])$.  Fix a point $p$ in the 
 interior of the  flat  sector $S$, and let $c_1$, $c_2$, $c_3$ be rays starting from $p$
   belonging  to  $\sigma(0)$, $\sigma(\pi)$  and  $\sigma(\pi/2)$ respectively.
Since $c_1$ and $c_3$ are contained in the  flat sector $S_1$, $\angle_p(c_1(\i), c_3(\i))=\pi/2$.
 Similarly  $\angle_p(c_2(\i), c_3(\i))=\pi/2$. Since the initial 
segment of each $c_i$ is contained in   the flat sector $S$, it   is   clear that the angle
 $\angle_p(c_1(\i), c_2(\i))=\pi$.  It follows that $c_1\cup c_2$ is a 
complete geodesic  in $X$  with endpoints  $\xi$    and  $\eta$.

\end{proof}

\b{remark}\label{a4coun}
\e{The conclusion  of Theorem \ref{pivisi}  does not hold  if $X$ is not a $CAT(0)$ $2$-complex.  For instance, 
 the universal covers of    nonpositively  curved $3$-dimensional graph manifolds (\cite{BS}, \cite{CK2})
%certain real analytic $4$-manifolds with nonpositive sectional curvature
 are counterexamples.}  % \e{(\cite{HS1}, \cite{AS})}.}
%There exist $CAT(0)$ spaces  $X$ and $\xi,\eta\in \p_T X$ with $d_T(\xi,\eta)=\pi$
\end{remark}

\section{Free Subgroups}\label{freesg}
    
  In  this section   we   use the results in Section   \ref{applic}
 to establish  a 
  sufficient condition (Theorem \ref{th4.11}) for the existence of free subgroups
  in a group acting isometrically %and cellularly 
on a $CAT(0)$   $2$-complex.

\subsection{Rank One Isometries} \label{rank1}

The reader is referred to the introduction   for 
  the definition of a hyperbolic isometry   $g$  and the notation $g(+\i)$, $g(-\i)$.

%is given in the introduction. 

\begin{Def}\label{rank1iso}
{A hyperbolic isometry $g$ of  a $CAT(0)$ space $X$ is  called  a \emph{rank  one isometry} %of $X$   
if   no  axis of $g$  bounds  a  flat half-plane in  $X$.}
\end{Def}

  We recall an  isometry  $g$ of a $CAT(0)$ space  $X$ induces a homeomorphism of   $\ol{X}$,
   which   we  still denote  by  $g$.

\begin{Th}{\emph{(W. Ballmann  \cite{B})}}\label{dyn}
{Let $X$ be a locally compact $CAT(0)$  space and 
  $g$      a rank one isometry of  $X$. 
   Given any neighborhoods
$U$ of $g(+\infty)$ and $V$ of $g(-\infty)$ in $\ol  X$,  there is an $n\ge 0$ 
such that  $g^k(\ol  X- V)\subset U$ and
$g^{-k}(\ol  X- U)\subset V$  whenever
$k\ge n$.}
\end{Th}

 Theorem \ref{dyn} implies  $g(+\infty)$ and $g(-\infty)$  are the only fixed points 
 of   a rank one isometry $g$ in $\ol  X$. 
   The theorem  also  has  the following   corollary.

\begin{Cor}\label{faf}
{Let $X$ be  a locally compact $CAT(0)$  space, 
  $G$     a  group of isometries of  $X$    and 
   $g\in G$    a rank one isometry.  %   with fixed points $g(+\infty)$,  $g(-\infty)$.
     Then one of the following holds:\newline
\emph{(i)}  $G$ has a fixed point in $\p_\i X$;\newline
 %$g(+\infty)$   or   $g(-\infty)$  is fixed by all elements of $G$; \newline
\emph{(ii)}  Some axis $c$ 
  of $g$  is  $G$-invariant; \newline 
\emph{(iii)}   $G$ contains  a free group of rank two.}
\end{Cor}

%\begin{Th}{\emph{(\cite{R})}}\label{kim}
%{ Let $X$ be a   locally   compact Hadamard space and $g$, $h$ two hyperbolic isometries of $X$. 
%If $d_T(\xi, \eta)>\pi$    whenever 
%   $\xi\in \{g(+\infty),g(-\infty)\}$ and $\eta\in \{h(+\infty),h(-\infty)\}$, 
%then the subgroup generated by
%$g$ and $h$     contains a free   group  of rank two.}
%\end{Th}

%It follows from  Theorem \ref{kim}  that   if    there are 
%two distinct Tits components
%$C_1$ and $C_2$ such that 
 %$\{g(+\infty), g(-\infty)\}\subset C_1$
  % and $\{h(+\infty), h(-\infty)\}\subset C_2$, 
 %then the group generated by $g$ and $h$ contains a free group of rank two.

\subsection{Ping-Pong Lemma} \label{ppl}

We recall the following well-known lemma.

\begin{Le}\label{pingpong}
{Let $G$ be a group acting on  a  set $X$, and  $g_1$, $g_2$   two  elements of $G$. 
If $X_1$ , $X_2$ are disjoint subsets of $X$ and for all $n\not=0$, $i\not=j$,
  $g_i^n(X_j)\subset X_i$, then the subgroup generated by $g_1$, $g_2$ is free of rank two.}

\end{Le} 

We will apply the Ping-Pong Lemma in the following setting. 
 Let $X$ be a $CAT(0)$ space and $g_1$, $g_2$ two hyperbolic isometries of $X$.
Then there are geodesics $c_1:R\ra X$, $c_2: R\ra X$ and  numbers
 $a, b>0$   with    $g_1(c_1(t))=c_1(t+a)$ and  $g_2(c_2(t))=c_2(t+b)$ for all $t\in R$.
Let $\pi_1: X\ra c_1(R)$ and $\pi_2: X\ra c_2(R)$ be orthogonal projections 
 onto the geodesics $c_1$ and $c_2$ respectively.   Set $X_1=\pi_1^{-1}(c_1((-\i,0]\cup [a, \i)))$ and
$X_2=\pi_2^{-1}(c_2((-\i,0]\cup [b, \i)))$.   The  following lemma is clear. 

\begin{Le}\label{proj}
{Let $X$, $g_1$, $g_2$, $X_1$ and $X_2$ be as above.  If $X_1\cap X_2=\phi$, then 
 the conditions in the Ping-Pong Lemma are satisfied. In particular, 
 $g_1$ and $g_2$ generate a free group of rank two.}
\end{Le}

\subsection{Free Subgroup Criterion} \label{fsc}

%Now we are ready to prove the main theorem of this section. 

A $CAT(0)$ $2$-complex is \e{piecewise Euclidean} if all its closed $2$-cells are isometric 
 to convex polygons in the Euclidean plane.
%The following is the main theorem of  Section  \ref{freesg}. 

\begin{Th}\label{criterion}
{Let $X$  be a   
  piecewise Euclidean  $CAT(0)$   $2$-complex  %  admitting a cocompact action by
  %cellular isometries  
  such that  the interior angles of   all the 
  $2$-cells are rational   multiples  of  $\pi$.    Suppose $g_1$ and $g_2$ are two hyperbolic isometries 
 of $X$ such that  for any $\xi\in \{g_1(+\i), g_1(-\i)\}$ and 
any $\eta\in \{g_2(+\i), g_2(-\i)\}$, there is a geodesic in $X$ with $\xi$   and 
  $\eta$   as endpoints.  
 Then  the group   generated by  $g_1$  and  $g_2$   contains  a free group 
 of rank two.}
\end{Th}

Let $c:I\ra X$ be a geodesic defined on an interval I. Then for any $p$ in the interior of $c(I)$, 
   the $d_p$ distance between  the two directions of $c$ 
at $p$ is  at least $\pi$.  Recall $d_p$ is a path metric defined on the link  $Link(X,p)$.   
 Notice in general  the two directions of $c$ 
at $p$  may have $d_p$  distance   strictly  larger than $\pi$. 

\begin{Def}\label{rgeo}
{A geodesic $c:I\ra X$ is an \e{R-geodesic} if for each point $p$ in 
the interior of $c(I)$,  the $d_p$  distance   between 
the two directions of $c$ 
at $p$  is $\pi$.}

\end{Def}

\begin{Le}\label{nonrgeo}
{With the assumptions of   Theorem \ref{criterion}.  If at least one of $g_1$, $g_2$ has 
 an axis that is not an  R-geodesic, then   the group   generated by  $g_1$  and  $g_2$   contains  a free group 
 of rank two.}
\end{Le}
\b{proof} Notice   if  a geodesic  bounds a flat half-plane then it 
is an   $R$-geodesic.
 The assumption in the lemma implies at least one of $g_1$, $g_2$ is a rank one isometry.
 % It is clear that 
  %if an axis of a hyperbolic isometry of a $CAT(0)$ 2-complex is not a $R$-geodesic, then
%this hyperbolic isometry is of rank one.  
Now the lemma follows   from   Corollary \ref{faf}.

\end{proof}

In light of  Lemma \ref{nonrgeo}, we will assume from now on that the axes of $g_1$ and $g_2$ are 
$R$-geodesics.
Since   $X$ admits a cocompact action by cellular isometries and  the interior angles of 
   all the $2$-cells are 
rational angles, there is a positive integer $m$ such that all the interior angles are integral  multiples
 of ${\pi}/{2m}$.

\b{Le}\label{perp}
{Given any two R-geodesic rays  $c_1:[0, \i)\ra X$, $c_2:[0,\i)\ra X$, there is a  number $a>0$ with the following 
 property:   if $q_i$ \e{($i=1,2$)} is a point in the interior of  $c_i$, and $p\in X$, $p\not=q_1, q_2$
 such that \newline
\e{(i)} $pq_1\cap pq_2=\{p\}$;\newline
\e{(ii)} $pq_1$ and $pq_2$ are R-geodesics;\newline
\e{(iii)} for each $i=1,2$, at least one of the angles $\angle_{q_i}(p, c_i(0))$,    
$\angle_{q_i}(p, c_i(\i))$  is $\pi/2$;\newline
 then $\angle_p(q_1, q_2)\ge a$.}
\end{Le}
\b{proof}
We first notice that for each $R$-geodesic $c$, there is a number $\alpha$, $0\le \alpha\le {\pi}/{4m}$  
  with the following property:
if $xy\subset c$ ($x,y\in c$) lies in some closed $2$-cell  $A$ of $X$, then  the angle between
 the segment $xy$ and  any direction at $x$ parallel to one of the edges of $A$ takes value in 
  the set   $\{\frac{k\pi}{2m}+\alpha: 0\le k< 2m\}\cup \{\frac{k\pi}{2m}-\alpha: 1\le k\le 2m\}$.  
 Let $\alpha_i$ ($i=1,2$) be such a number corresponding to $c_i$.   Since $pq_i$ is also 
an  $R$-geodesic and at least one of the angles $\angle_{q_i}(p, c_i(0))$,    
$\angle_{q_i}(p, c_i(\i))$  is $\pi/2$,   the number $\alpha_i$ also corresponds to 
  $pq_i$.  The fact  $pq_1\cap pq_2=\{p\}$  implies $\angle_p(q_1, q_2)\not=0$. 
Now it   is   easy to see that 
 $\angle_p(q_1, q_2)\ge |\alpha_1-\alpha_2|$  when $\alpha_1\not=\alpha_2$;
$\angle_p(q_1, q_2)\ge \min\{2\alpha_1, \frac{\pi}{2m}-2\alpha_1\}$   when 
 $0<\alpha_1=\alpha_2<{\pi}/{4m}$;   
 and $\angle_p(q_1, q_2)\ge \pi/{2m}$   when 
$\alpha_1=\alpha_2=0$   or ${\pi}/{4m}$.   

\end{proof}

\b{Le}\label{smallangle}
{Let $\xi_1, \xi_2\in \p_\i X$    and 
$c_1, c_2: [0, \i)\ra X$ \e{($i=1,2$)}  rays belonging to $\xi_1$  and $\xi_2$  respectively.
    Suppose 
there is a geodesic $c$  in $X$  with  $\xi_1$   and  $\xi_2$ as endpoints.  
Then for $i\not=j$,  $\angle_{c_i(t)}(c_i(0),  c_j(t))\ra 0$ as $t\ra \i$.}
\end{Le}
\b{proof}
Since $c$ is a  geodesic  in $X$  with  $\xi_1$    and  $\xi_2$  as endpoints,  there is a number $\epsilon>0$ such that
 $c_i\subset N_\epsilon(c)$ ($i=1,2$).   The convexity of distance function implies 
     for large enough $t$, there is a 
point $p_t\in c_1(t)c_2(t)$ with 
$d(c(0), p_t)\le \epsilon$.  Now the lemma follows by considering the triangle
 $\Delta(c_i(t)c_i(0)p_t)$. 

\end{proof}

 Let $g_1$, $g_2$ be as in Theorem \ref{criterion} and $c_1$, $c_2$  be 
axes of $g_1$, $g_2$ respectively. By Lemma \ref{nonrgeo}, we may  assume 
$c_1$, $c_2$   are $R$-geodesics. 
Let $\pi_i:X\ra c_i$ ($i=1,2$) be the orthogonal 
projection onto $c_i$.
 If there is some  $T>0$ such that  
$$\pi_1^{-1}(c_1((-\i, -T]\cup [T, \i)))\cap \pi_2^{-1}(c_2((-\i, -T]\cup [T, \i)))=\phi,$$ 
  then for large enough $n$, $g_1^n$ and $g_2^n$ satisfy the condition in Lemma \ref{proj}
  and  therefore generate a free group of rank two.
We shall prove  there is some $T>0$   with  $\pi_1^{-1}(c_1([T, \i)))\cap \pi_2^{-1}(c_2([T, \i)))=\phi$, 
  the other three cases are similar. 

\b{Le}\label{emptyi}
{Let $\xi_1, \xi_2\in \p_\i X$    and 
   $c_1, c_2: [0, \i)\ra X$ \e{($i=1,2$)}  rays belonging to $\xi_1$  and $\xi_2$  respectively.
 Suppose 
  $c_1$ and $c_2$ are $R$-geodesics and 
there is a geodesic $c$ in $X$ with   $\xi_1$   and $\xi_2$  as endpoints.   Then  
$\pi_1^{-1}(c_1([T, \i)))\cap \pi_2^{-1}(c_2([T, \i)))=\phi$  
for large enough $T>0$, 
where 
 $\pi_i: X\ra c_i$ is the orthogonal projection onto $c_i$.}
\end{Le}
\b{proof}
Since  $c_1$ and $c_2$ are $R$-geodesics, there is a  number $a>0$ with the 
 property  stated in Lemma \ref{perp}.  By Lemma \ref{smallangle}, 
 there is  $T>0$   with  $\angle_{c_1(t)}(c_1(0),  c_2(t))<a/3$  and 
$\angle_{c_2(t)}(c_2(0),  c_1(t))<a/3$  for all $t\ge T$.

Suppose  there  exists  $t>T$ with
$\pi_1^{-1}(c_1([t, \i)))\cap \pi_2^{-1}(c_2([t, \i)))\not=\phi$. We will derive a contradiction from this. 
  Pick
$p\in \pi_1^{-1}(c_1([t, \i)))\cap \pi_2^{-1}(c_2([t, \i)))$ and set $q_i=\pi_i(p)$ ($i=1,2$). 
Let $\sigma_2: [0, b_2]\ra X$ be the geodesic from $p$ to $q_2$.
   Consider  ${\pi_1}_{|{\sigma_2}}: \sigma_2\ra c_1$ and 
let $t_2=\max\{t\in [0, b_2]:\pi_1(\sigma_2(t))=q_1\}$.  Then the geodesic $q_1\sigma_2(t_2)$
  is an   $R$-geodesic and at least one of the angles 
 $\angle_{q_1}(c_1(0), \sigma_2(t_2))$, $\angle_{q_1}(c_1(\i), \sigma_2(t_2))$ is $\pi/2$.
   By replacing $p$ with $\sigma_2(t_2)$, we may assume  $pq_1$ is an   $R$-geodesic and 
 at least one of the angles $\angle_{q_1}(c_1(0), p)$, $\angle_{q_1}(c_1(\i), p)$ is $\pi/2$.
  Similarly we may assume 
$pq_2$ is an   $R$-geodesic and 
 at least one of the angles $\angle_{q_2}(c_2(0), p)$, $\angle_{q_2}(c_2(\i), p)$ is $\pi/2$.
In general  $pq_1$ and $pq_2$ share an initial segment: $pq_1\cap pq_2=pp'$. By replacing 
  $p$ with $p'$ if necessary we may assume $pq_1\cap pq_2=\{p\}$.  Now the conditions of 
Lemma \ref{perp} are satisfied and  so $\angle_p(q_1, q_2)\ge a$.

  Consider the triangle $\Delta(pq_1q_2)$.  Since $q_1=\pi_1(p)$,   we   have  
$\angle_{q_1}(p, c_1(0))\ge \pi/2$. It follows that 
$\angle_{q_1}(p, q_2)\ge \angle_{q_1}(p, c_1(0))-\angle_{q_1}(q_2, c_1(0))\ge \pi/2-a/3$.
  Similarly $\angle_{q_2}(p, q_1)\ge  \pi/2-a/3$.
  Now $\angle_{q_1}(p, q_2)+ \angle_{q_2}(p, q_1)+\angle_p(q_1, q_2)\ge 2(\pi/2-a/3)+a>\pi$,
   a contradiction.

\end{proof}

The proof of Theorem \ref{criterion}
 is now complete.

The following result is an immediate consequence of Theorems \ref{criterion}
   and \ref{pivisi}.

\begin{Th}\label{th4.11}
{Let   $X$ be a    piecewise Euclidean $CAT(0)$   $2$-complex    
%that   admits a   cocompact   action by cellular isometries 
   so that  the  interior 
angles of   all the $2$-cells  are rational multiples of  $\pi$, and 
   $g_1$,  $g_2$  two hyperbolic isometries 
 of $X$.     If   $g_1(+\i)$,    $g_1(-\i)$,   $g_2(+\i)$,   $g_2(-\i)$  are not terminal points   and 
  $d_T(\xi, \eta)\ge \pi$ for 
 any $\xi\in \{g_1(+\i), g_1(-\i)\}$    and any 
   $\eta\in \{g_2(+\i), g_2(-\i)\}$,   
   then  the group   generated by  $g_1$  and  $g_2$   contains  a free group 
 of rank two.}
\end{Th}

\begin{Cor}\label{tl}
{Let   $X$ be a    piecewise Euclidean   $CAT(0)$   $2$-complex   such that the  interior 
angles of    all the $2$-cells  are rational   multiples  of  $\pi$,
and $G$ a group acting   on  $X$ properly  and     
    cocompactly   by     cellular isometries.   Suppose 
   $g_1, g_2\in G$   are   two hyperbolic isometries 
   such that   $d_T(\xi, \eta)\ge \pi$ for 
 any $\xi\in \{g_1(+\i), g_1(-\i)\}$    and any 
   $\eta\in \{g_2(+\i), g_2(-\i)\}$. 
 Then  the    subgroup   generated by  $g_1$  and  $g_2$   contains  a free group 
 of rank two.}
\end{Cor}

\begin{proof}
Notice in this case for any hyperbolic isometry   $g\in G$   either 
 $g$ is a rank one isometry or $g(+\i)$,  $g(+\i)$  are not terminal points.  
 The corollary follows from   Theorem  \ref{th4.11}    and   Corollary \ref{faf}.

\end{proof}

\section{Quasi-isometries Between CAT(0) 2-complexes}  \label{sec4}

   In this section  we  study how the Tits boundary behaves under quasi-isometries.   
Throughout this section, 
   $CAT(0)$ 2-complexes  %locally compact $CAT(0)$     piecewise Riemannian  2-complexes 
 %that admit   cocompact isometric actions  
  are as defined in  Section \ref{2complex}.

  For any    $CAT(0)$  $2$-complex  $X$,   
%locally compact $CAT(0)$     piecewise Riemannian  2-complex that admits a  cocompact   isometric 
   %group action. 
set 
$ Core(\partial_T  X)=\cup c$
 where $c$   varies  over all the topological circles in $\partial_T X$.  Let $d_c$ be the induced path metric  of $d_T$  on 
 $Core(\partial_T X)$.    Then $d_c\ge d_T$  always holds.    %Since   $\partial_T X$ is a $CAT(1)$ space  
 %and has geometric dimension  at most  1,
  %  any  closed   metric ball   in    $\partial_T X$   with radius  $\frac{\pi}{4}$ is an    $R$-tree (Corollary \ref{c1}).  
 %It follows that 
   Recall by  Corollary \ref{c1}  $\ol B(\xi,r)\subset \p_T X$ is  an $R$-tree 
for any $\xi\in  \p_T X$ and any  $r$  with $0<r<\pi/2$.   It follows that   
      for $\xi, \eta\in Core(\partial_T  X)$,    
%$d_c(\xi, \eta)=d_T(\xi, \eta)$  when  $d_c(\xi, \eta)< \pi/2$.  It is also not hard to see that
$d_c(\xi,\eta)<\i$  if  and only if  $\xi$, $\eta$ lie in the same path component of $Core(\p_T X)$
   and in this case there is  a minimal Tits geodesic  contained in  $Core(\p_T X)$   and 
connecting   $\xi$ and  $\eta$. In particular,   
   $d_c(\xi,\eta)=d_T(\xi,\eta)$ if $d_c(\xi,\eta)<\i$.   

%$\xi$, $\eta$ lie in the same path component of $Core(\p_T X)$, and $d_c(\xi,\eta)=\i$ otherwise. 
 
Below is the main result of this section.

\begin{Th} \label{t17}
{Let $X_1$ and $X_2$ be two    $CAT(0)$  $2$-complexes. If $X_1$ and $X_2$  
are $(L,A)$ quasi-isometric, then $Core(\partial_T X_1)$ and $Core(\partial_T X_2)$  are  $ L^2$-bi-Lipschitz with respect to the metric $d_c$. }
\end{Th}

\subsection{Tits Limit Set}  \label{tlimit}

Let $X$ be  a   $CAT(0)$ $2$-complex.  For any subset $I\subset \partial_T X$ and 
$x\in X$, the geodesic cone over $I$ 
  with vertex $x$  is  $C_x(I)=\bigcup_{\xi\in I}x\xi$.   Let 
${\mathcal C}(\partial_T X)$   be  the set of topological circles in $\partial_T X$.
 %For  $S\in {\mathcal C}(\partial_T X)$    and   $x\in X$, the geodesic cone over $S$ 
 % with vertex $x$  is  $C_x(S)=\bigcup_{\xi\in S}x\xi$. 
   By Corollary \ref{r1}   each  $S\in {\mathcal C}(\partial_T X)$ 
 is a simple closed geodesic in $\partial_T X$.  It is not hard to see that
 $C_x(S)$  is  a quasi-flat.  
  The   main ingredient   in   the proof of  Theorem \ref{t17}
is a result of B. Kleiner concerning top dimensional quasi-flats in $CAT(0)$   spaces. A 
special case of his result is as follows:

 \begin{Th}\label{t18}
\e{(Kleiner \cite{K2})}
{Let  $X$ be  a  $CAT(0)$   $2$-complex   and 
 $Q\subset X$  a 
quasi-flat.   Then there is a unique  $S\in {\mathcal C}(\partial_T X)$   such that  
    for any $x\in X$,   
  $$\lim_{r\ra \i}\frac{d_H(Q\cap B(x, r), C_x(S)\cap B(x, r))}{r}=0.$$  }

\end{Th}

For convenience, we
 introduce the following definition.
 
\begin{Def}\label{d8}
{Let $B\subset X$ be a subset of a $CAT(0)$  space $X$. A point $\xi\in \partial_T X$ is  a \emph{Tits limit point}  of $B$ 
 if there is a sequence of points $b_i\in B$,    $i=1,2,\cdots, $ such that  $d(b_i,x)\ra \i$  and 
 $\lim_{i\rightarrow \infty}{\frac{d(b_i,x\xi)}{d(b_i,x)}}=0$ 
where  $x\in X$ is  a fixed point.    We also say  $\{b_i\}_{i=1}^\i$ \e{Tits converges}  to $\xi$. 
   The   Tits limit   set $L_T(B)$  of  $B$  is  the set 
of  Tits limit points of $B$.     } 
\end{Def}

 Notice the above definition does not depend on the choice of $x$. 
  It is clear that $L_T(B)$  is closed in the Tits metric,  and if 
  $d_H(B_1, B_2)<\i$ for  $B_1,B_2\subset X$ then  $L_T(B_1)=L_T(B_2)$.
  A  Tits limit point  is a limit point in the usual sense, but 
 a limit point does not have to   be a Tits limit point.  
  Limit points are 
defined in terms of the cone topology, while  Tits limit points  are      related to the Tits metric.
Using this terminology, Kleiner's theorem   in particular implies the following:  
$L_T(Q)\in  {\mathcal C}(\partial_T X)$  
for  any quasi-flat  
$Q\subset X$.    We   also   have $L_T(E)=C$ for the quasi-flat in  Proposition  \ref{c3}.

\begin{Le}\label{l32}
{Let $S\in {\mathcal C}(\partial_T X)$, 
      $p\in X$  and  $x_i\in  X$ \e{($i=1,2,\cdots$)} be  a sequence of 
points with 
$lim_{i\rightarrow \infty}{d(x_i, p)}=\infty$ and
 $lim_{i\rightarrow \infty}{\frac{d(x_i,C_p(S))}{d(x_i,p)}}=0$.   Then  some 
subsequence of $\{x_i\}$ Tits converges to a point $\xi\in S$.}
\end{Le}
 
\begin{proof}     Let $y_i\in C_p(S) $ with $d(x_i,y_i)=d(x_i,C_p(S))$. 
  Then $y_i\in p\xi_i$  for some $\xi_i\in S$. 
Since $S$ is a simple closed geodesic in $\p_T X$, a subsequence $\{\xi_{i_j}\}$   of   $\{\xi_i\}$ 
  converges to some $\xi\in S$ in the Tits metric.    
%Let   $\sigma$  be a ray    of    the foliation representing $\xi$.  
 Then $\lim_{{j}\rightarrow 
\infty}{\frac { d(y_{i_j}, \; p\xi)}{d(y_{i_j},\;  p)}}=0$.   Now   triangle inequality implies  $\lim_{{j}\rightarrow \infty}{\frac {d(x_{i_j}, \; p\xi  )}{d(x_{i_j},  \; p)}}=0$ and $\{x_{i_j}\}$ Tits converges to $\xi$. 

\end{proof}

\begin{Def}\label{d9}
{Let   $L\ge 1$, $A>0$.
A metric space $M$ is \e{$(L,A)$ quasi-connected at infinity}  if there is some  point $a\in M$  and
 some      $r_0>0$  with the following property:  
 for any two points 
 $x,y\in M$  with $d(x,a), d(y,a) > r_0 $  there is  a sequence of points $x=x_0, x_1,\cdots, 
x_k=y$      so that:  \newline
(i) $d(x_i,x_{i+1})\le A$ for $i=0, 1,   \cdots, k-1$;\newline
(ii)  $d(x_i, a)\ge {\min\{d(x,a), d(y,a)\}}/L$  for all $i$.\newline
 A metric space $M$ is  {quasi-connected at infinity} if it is 
   $(L,A)$ quasi-connected at infinity for some $L\ge 1$, $A>0$.}
\end{Def}

Notice    %quasi-connectedness at infinity   is a quasi-isometry invariant,
    if two metric spaces $M_1$  and $M_2$ are quasi-isometric,  then
 $M_1$  is quasi-connected at  infinity if and only  if $M_2$ is.

\begin{Le}\label{l33}
{Let  $X$ be  a $CAT(0)$ $2$-complex and  $Q\subset X$ a quasi-flat. 
 If  $B\subset Q$ 
 is      quasi-connected at infinity,   
  then  $L_T(B)\subset  L_T(Q)$ is   path-connected.  }
\end{Le}
\begin{proof} Suppose   $B$ 
 is     $(L,A)$ quasi-connected at infinity.
 Set $S=L_T(Q)$.  
$L_T(B)\subset S$ is clear since $B\subset Q$.  
 Let $\xi\not=\eta\in L_T(B)$. 
 %Since $L_T(B)\subset S$ there are two rays  $\alpha$  and  $\beta$  from the foliation
%that represent  $\xi$  and  $\eta$  respectively.  
 There are two sequences  $\{a_i\}\subset B$,  $\{b_i\}\subset B$    with      $d(a_i,a)\ra \i$,  
$\frac{d(a_i, a\xi)}{d(a_i,a)}\rightarrow 0$    and  $d(b_i,a)\ra \i$,  
 $\frac{d(b_i, a\eta)}{d(b_i,a)}\rightarrow 0$ 
as $i\rightarrow \infty$,  where $a$ is the base point of $B$ in
   the definition of quasi-connectedness.     The quasi-connectedness of $B$ implies  for each $i$  there 
is a sequence  $a_i=x_i^1, x_i^2, \cdots, x_i^{k_i}=b_i$,   $x^j_i\in B$ with 
$d(x_i^j, a)\ge {\min\{d(a_i,a), d(b_i,a)\}}/L$    and 
$d(x_i^j,x_i^{j+1})\le A$.          Let $\tilde x_i^j\in C_a(S)$
  with  $d(x_i^j, \tilde x_i^j)=d(x_i^j, C_a(S))$.
 %be the  point in $\ol E$ closest  to $x_i^j$. 
   We also let $\tilde x_i^0\in  a\xi $,  $\tilde x_i^{k_i+1}\in a\eta $
  with $d(a_i, \tilde x_i^0)=d(a_i, a\xi)$  and 
$d(b_i, \tilde x_i^{k_i+1})=d(b_i,  a\eta)$.  
  Then  $\tilde x_i^j\in a\xi_i^j$   for some  
 $\xi_i^j\in S$.     
   We may  choose      $\xi_i^0=\xi$,  $\xi_i^{k_i+1}=\eta$. 
 The facts   $d(x_i^j,x_i^{j+1})\le A$   and  
${d_H(Q\cap B(a, r), C_a(S)\cap B(a,r))}/r\rightarrow 0$ as $r\rightarrow \infty$
       imply  $\frac{d(\tilde x_i^j,\tilde x_i^{j+1})}{d(\tilde x_i^j,a)}\ra 0$  as $i\ra \i$.
  Since $S$ is  compact  in $\p_T X$ it is not hard to see that 
there are positive numbers $\epsilon_i$ with $\epsilon_i\rightarrow 0$ as $i\rightarrow \infty$ and
$d_T(\xi_i^j,\xi_i^{j+1})<\epsilon_i$.

 $\xi $ and $\eta $ divide $S$ into two closed intervals $I_1$ and $I_2$.  After possibly passing to a subsequence and relabeling 
 $I_1$ and $I_2$ we have  $N_{\epsilon_i}(\{\xi_i^0, \cdots, \xi_i^{k_i+1}\})\supset I_1$.  Now it  follows 
 from the definition of
 Tits limit points that   $I_1\subset L_T(B)$ and 
$\xi$, $\eta$  can be connected by a path in $L_T(B)$.  

 \end{proof} 

\subsection{Induced Map Between Sets of Branch Points}  \label{indmbbp}

  %For  any  $CAT(0)$ $2$-complex  $X$,   let 
%$P(\partial_T X)$ be the set of topological circles in $\partial_T X$.   %Recall   given any  
 %$S\in   P(\partial_T X)$, Corollary \ref{c3}   provides a quasi-flat $E$ with $L_T(E)=S$.
 Let  $f:X_1\rightarrow X_2$   be  a  $(L,A)$  quasi-isometry between $CAT(0)$ 2-complexes. 
Then   $f$ induces a 
bijection $g: {\mathcal C}(\partial_T X_1)\rightarrow {\mathcal C}(\partial_T X_2)$  as follows. 
   For  each  circle  $S\in {\mathcal C}(\partial_T X_1)$,  
  define   $g(S)=L_T(f(C_x(S)))$ for any $x\in X_1$.   It  is clear $g(S)$  does not depend on $x$. 

   To prove  Theorem  \ref{t17},  we shall  first show that the  bijection $g$ preserves the 
intersection pattern of circles. This implies there is a bijection between the set of branch 
points of $Core(\partial_T X_1)$ 
and that of 
$Core(\partial_T X_2)$. We then show that this bijection between
     branch points is actually induced by the quasi-isometry $f$,
   so it is a local   bi-Lipschitz map (see Definition \ref{d10}).  Finally we 
    extend this map to $Core(\partial_T X_1)$.

First we look at the intersection of two circles in $\partial_T X$.
%For any topological space $Y$, $\pi_0(Y)$ denotes the set of path components of  $Y$.  

\begin{Le}\label{l31}
{Let  $X$  be  a  $CAT(0)$ $2$-complex   and 
 $S_1,S_2\in {\mathcal C}(\partial _TX)$. Then $S_1\cap S_2\subset  S_i$ is a disjoint union of finitely many points and 
finitely many closed intervals.}
\end{Le}
\begin{proof}
 Write $S_1=\cup_i{U_i}$ as   a  finite  union of closed intervals each of which has length $< \pi$. Similarly  
$S_2=\cup_j{V_j}$. Then  $S_1\cap S_2=\cup_{i,j}(U_i\cap V_j)$.   Since $\partial _T X$ is a $CAT(1)$ space, and $U_i$ and 
$V_j$ are closed intervals of length $<\pi$, $U_i\cap V_j$  is empty, a single point or a closed 
interval.  Thus $S_1\cap S_2=\cup_{i,j}(U_i\cap V_j)\subset  S_i$ is a union of finitely many points and finitely many closed 
intervals.   

 \end{proof}

 For any topological space $Y$, $\pi_0(Y)$ denotes the set of path components of  $Y$. 
%We  fix  a quasi-inverse $f^{-1}:X_2\ra  X_1$  of $f$. 

 \begin{Le}\label{l34}
  {Let    $S_1,S_2\in {\mathcal C}(\partial _T X_1)$.  %\in {\mathcal C}(\partial_T X_1)$.  
Then there   is  a unique bijective map
  $k:\pi_0(S_1\cap S_2)\ra  \pi_0(g(S_1)\cap g(S_2))$  such that  $L_T(f(C_x(I)))=k(I)$
for  all $I\in \pi_0(S_1\cap S_2)$,  where  $x\in X_1$.} 
  %    and   $L_T(f^{-1}(C_y(B)))\subset  k^{-1}(B)$  for 
%$B\in \pi_0(g(S_1)\cap g(S_2))$,  where  $x\in X_1$ and $y\in  X_2$.}

\end{Le}
 
\begin{proof}  
 Let $x\in X_1$ and    $I\in   \pi_0(S_1\cap S_2)$.    $C_x(I)$  is  quasi-connected at  infinity 
 since  by Lemma  \ref{l31}  $I$ is a single point or a closed interval.     
  $f$  is  a quasi-isometry  implies  
$f(C_x(I))$  is    also  quasi-connected at infinity.   Since 
 $C_x(I)\subset C_x(S_i)$ ($i=1,2$) we  have  $f(C_x(I))\subset f(C_x(S_i))$  and 
 $L_T(f(C_x(I)))\subset L_T(f(C_x(S_i)))=g(S_i)$.
It follows from Lemma \ref{l33}  that    $L_T(f(C_x(I)))\subset g(S_1)\cap g(S_2)$ is path  connected.
  Let $k(I)$ be the component of  $g(S_1)\cap g(S_2)$  that  contains  $L_T(f(C_x(I)))$.  
  By  considering    a  quasi-inverse      of $f$  and  using  Lemma  \ref{l32}, 
we see $L_T(f(C_x(I)))=k(I)$
   and  that 
 $k$  is  bijective.    $k$ is clearly     unique.  

  %Since $f:X_1\rightarrow X_2$ is a quasi-isometry, by   considering 
  %a    quasi-inverse of $f$,  we see that  $L_T(f(B_1))=A'$.    

\end{proof}

\begin{Le}\label{l35}
{%With the notation of  Lemma  \ref{l34}.
Let $I\in \pi_0(S_1\cap S_2)$.  Then $I$ is a  single point   component 
  if and only if 
  $k(I)$ is. }
 \end{Le}

\begin{proof} Suppose the lemma is false.  By possibly replacing $f$ with a quasi-inverse, we may assume 
  $g(I)=\{\eta\}$ ($\eta\in \p_T X_2$)  for some 
  nontrivial closed interval component $I$ of $ S_1\cap S_2$.
 %Fix some  $x\in  X_1$. 
  By Lemma  \ref{l34}  
$L_T(f(C_x(I)))=\{\eta\}$,  where $x\in  X_1$ is a fixed  base point.  
  Then  there are   $\xi_1\not=\xi_2\in   I$   with 
$L_T(f(x\xi_1))=L_T(f(x\xi_2))=\{\eta\}$.
    Let  $x_i\in x\xi_1$  and  $y_i\in x\xi_2$
  ($ i\geq 1$)  with  $d(x_i,x)=d(y_i,x) =i$.   Then   both $\{f(x_i)\}$ and $\{f(y_i)\}$ Tits converge to 
$\eta$. 
Let $x_i',  y_i'\in f(x)\eta$ ($i\ge 1$)    with $d(f(x_i), x_i')=d(f(x_i), f(x)\eta)$  and   
   $d(f(y_i), y_i')=d(f(y_i),f(x)\eta)$.   
 Then for  large enough $i$, there is 
some $j_i$ such that $d(x_i',y'_{j_i})\leq L+A$ ($f$ is a $(L,A)$ quasi-isometry).
  Then $\lim_{i\rightarrow \infty}{\frac{d(f(x_i),  f(y_{j_i}))}{d(f(x_i),  f(x))}}=0$. 
%where $p\in X_1$ is a fixed basepoint. 
 Since $f$ is a quasi-isometry, we have  $\lim_{i\rightarrow \infty}{\frac{d(x_i,  y_{j_i})}{d(x_i,x)}}=0$. This 
contradicts to the assumption that $x_i\in x\xi_1$, $y_{j_i}\in x\xi_2$ and 
 $\xi_1\not=\xi_2$.  

 \end{proof}

 %Lemma   \ref{l35} implies 
  % $k$  induces  a map from  the set of single point components of $S_1\cap S_2$   is Tits induced by  $f$.

 \begin{Le}\label{l36}
 { Let $I$ be a  nontrivial  interval component of  $S_1\cap S_2$, $\xi_1, \xi_2$  the two 
 endpoints of $I$, and $\eta_1, \eta_2$ the two endpoints of $k(I)$. 
 Then  for any  $x\in X_1$ either $L_T(f(x\xi_1))=\{\eta_1\}$, 
$L_T(f(x\xi_2))=\{\eta_2\}$, or
$L_T(f(x\xi_1))=\{\eta_2\}$, $L_T(f(x\xi_2))=\{\eta_1\}$.}
 \end{Le}
\begin{proof}  We   claim    $L_T(f(x\xi_i))$ does not contain any interior point of $k(I)$.
Suppose     at least one of   $L_T(f(x\xi_1))$,  $L_T(f(x\xi_2))$, say 
 $L_T(f(x\xi_1))$    does  contain     some interior point $\eta $ of 
$k(I)$. Let
 $d_0=d_T(\eta, g(S_1)-k(I))$. 
Choose $\xi\in S_1-I$ with $d_T(\xi_1, \xi)$ sufficiently small.     
Since $\eta\in L_T(f(x\xi_1))$, there is a sequence of points $x_i\in x\xi_1$, such that $\{f(x_i)\}$ Tits converges 
to $\eta$. Let $y_i\in x\xi$   with  $d(x, y_i)=d(x,x_i)$.  By passing to a subsequence, we may assume that 
$\{f(y_i)\}$ Tits converges to a point $\eta'\in g(S_1)$.   Since $f$ is a  quasi-isometry 
  we see  $d_T(\eta',\eta)\le {d_0}/{2}$ if $d_T(\xi_1, 
\xi)$ is sufficiently   small. 
The choice of $d_0$ then implies   $\eta'\in k(I)$. 
  It  follows  that  $\xi\in L_T(f^{-1}(C_y(k(I))))$, where $y\in X_2$  and 
$f^{-1}: X_2\rightarrow X_1$    is  
a  quasi-inverse
of $f$.   This  is  a contradiction since $L_T(f^{-1}(C_y(k(I))))=I$.

 By  Lemma \ref{l33}  
   $L_T(f(x\xi_i))$ is          path connected.   
So  either 
    $L_T(f(x\xi_i))=\{\eta_1\}$   or  $L_T(f(x\xi_i))=\{\eta_2\}$. 
Similarly  
$L_T(f^{-1}(y\eta_i))=\{\xi_1\}$ or $\{\xi_2\}$.           Now the lemma follows easily.

  %So  either 
   % $L_T(f(x\xi_i))=\{\eta_1\}$   or  $L_T(f(x\xi_i))=\{\eta_2\}$.   But 
 %$L_T(f(x\xi_1))=\{\eta_1\}=L_T(f(x\xi_2))$ can not occur  since it would imply 
%$L_T(f^{-1}(y\eta_1))=\{\xi_1, \xi_2\}$ is disconnected, where $f^{-1}$ is a quasi-inverse of $f$.    Similarly, 
 %$L_T(f(x\xi_1))=\{\eta_2\}=L_T(f(x\xi_2))$  does not happen and the lemma follows.   

  \end{proof}

\begin{Def}\label{tinduced}
Let  $f:X_1\rightarrow X_2$  be  a  quasi-isometry between two  $CAT(0)$ $2$-complexes
  and $A_1\subset \p_T X_1$, $A_2\subset \p_T X_2$.  A map $h: A_1\ra A_2$ is 
  \e{Tits induced} by  $f$ if  it satisfies the following property: 
    for  any   sequence $\{x_i\}\subset X_1$ Tits converging     to $\xi\in A_1$, 
    $\{f(x_i)\}$ Tits converges  
     to $h(\xi)$.

\end{Def}

  Define  $B_1\subset Core(\p_T X_1)$ as follows:
 $$B_1=\{\xi\in \p_T X_1: \xi \;\text{is an endpoint of some}\; 
I\in \pi_0(S_1\cap S_2), S_1, S_2\in {\mathcal C}(\p_T X_1)\}.$$
 Here $\xi$  is an endpoint of $I$   when  $I=\{\xi\}$. 
Similarly we define $B_2\subset Core(\p_T X_2)$.  Lemma \ref{l35} and
Lemma \ref{l36} imply the following  proposition:

\begin{Prop}\label{indbbp}
{Let $f:X_1 \rightarrow X_2$  be  a 
 quasi-isometry   between   two  $CAT(0)$ $2$-complexes, and $B_1$, $B_2$ as above.
  Then there is  a bijective map $h:B_1\ra  B_2$   with 
 $h(S\cap B_1)=g(S)\cap B_2$  for any 
  $S\in  {\mathcal C}(\p_T X_1)$,  such  that $h$ is Tits induced by $f$ and 
 $h^{-1}$  is   Tits  induced by a quasi-inverse of $f$. }

\end{Prop}

   By  using  Lemma \ref{l33}
  and  considering  a  quasi-inverse of  $f$ it is  not hard to show the following  lemma.
 %The following corollary follows easily from Lemma \ref{l36}.
 
%\begin{Cor}\label{c4}
 %{The bijection in Lemma \ref{l34}  induces a bijective map between the set of endpoints of closed interval components of 
%$S_1\cap S_2$ and that of $g(S_1)\cap g(S_2)$, and this bijective map is induced by the quasi-isometry $f:X_1\rightarrow 
%X_2$. }
%\end{Cor}
 
%It follows from this corollary and the remark after Lemma \ref{l35}  that the quasi-isometry $f:X_1 \rightarrow X_2$ induces a 
%bijective map $h$ from the set of branch points $B_1$ in $Core(\partial_T X_1)$ to
%the set of branch points $B_2$ in $Core(\partial_T X_2)$. Furthermore, for any topological circle $S\subset \partial_T 
%X_1$, $h(S\cap B_1)=g(S)\cap B_2$.

\begin{Le}\label{l37}
{Let $S\in {\mathcal C}(\p_T X_1)$.   Then   $h_{|S\cap B_1}: S\cap  
B_1\rightarrow g(S)\cap B_2$ preserves the order of points, that is,   
   if $a_1, a_2,a_3,a_4\in S\cap B_1 $ are four 
points in cyclic  order    on   $S$, then $h(a_1),
h(a_2),h(a_3),h(a_4)$ are in cyclic  order   on   $g(S)$. }
\end{Le}
%\begin{proof}     We first  assume 
%$L_T(f(B))$   contains  $h(a_3)$.        Then there is a sequence $x_i\in B$ so that  $\{f(x_i)\}$ Tits converges to  $h(a_3)$.  
%After passing to a subsequence we may  assume $\{x_i\}$ Tits converges to some $\xi\in I$.   Let $c'\subset E'$ be a ray representing 
 %$h(a_3)$.  Since  $\{x_i\}$ Tits converges to     $\xi\in I$ and  $\{f(x_i)\}$ Tits converges to  $h(a_3)$,  we see $\xi\in L_T(f^{-1}(c'))$
  %where $f^{-1}$ is a quasi-inverse of $f$.   Here comes the  contradiction since $ L_T(f^{-1}(c'))=\{a_3\}$,  $a_3\notin I$ 
  %and $\xi\in I$.    Similarly we reach a contradiction when $L_T(f(B))$   contains  $h(a_4)$.  

 %\end{proof} 

\subsection{Bi-Lipschitz Map  Between the Cores}  \label{localbl}

Let  $f:X_1\rightarrow X_2$  be  a  $(L,A)$  quasi-isometry.
In this  section we shall extend the map $h$ in Proposition \ref{indbbp}
 to  a bi-Lipschitz map from $Core(\p_T X_1)$ to $Core(\p_T X_2)$. 
 We first show $h$ is a local bi-Lipschitz map.  

\begin{Def} \label{d10}
{ Let $L_0\ge 1$.   
 A map $h:Y_1\rightarrow Y_2$ between two metric 
spaces is    a   \e{ local   $L_0$-bi-Lipschitz map} if  
 there is some $\epsilon>0$  such that 
$$1/{L_0}\cdot d(a,b)\leq d(h(a),h(b))\leq L_0\cdot d(a,b)$$ 
for all $a,b\in Y_1$ with $d(a,b)\leq \epsilon$.}
\end{Def}

The following lemma is easy to prove. 

\begin{Le}\label{l38}
{ Let $Y_1$,  $Y_2$ be two metric spaces,   $A_1\subset Y_1$,     and    
  $h:A_1\rightarrow Y_2$  a   local    $L_0$-bi-Lipschitz  map for some $L_0\ge 1$. 
   If $Y_2$ is  complete, 
 then $h$   uniquely extends to a    
local   $L_0$-bi-Lipschitz  map  $\bar h: \bar A_1\rightarrow    Y_2$.}
\end{Le}

\begin{Le}\label{l39}
{For any  $\lambda >1$,        the map  $h:B_1\rightarrow B_2$ is a   local  $\lambda L^2$-bi-Lipschitz    map   with respect to the 
  Tits metric $d_T$.  }
\end{Le}
\begin{proof}   Given  $\lambda >1$,    fix    some $\mu>0$ so that  if   $\sin t\le  \mu$,     $ 0<t<{\pi}/{2}$ then  
 ${t}/{\lambda} \le \sin t \le t$.  
 Set  $\epsilon= {\mu}/{L^2}$.      Let $\xi, \eta\in B_1$ with $d_T(\xi,\eta)\leq \epsilon$.  
  Fix  a point $p\in X_1$ and let $x_i=\gamma_{p\xi}(i)$,  $y_i=\gamma_{p\eta}(i)$, $i\ge 1$.
Then $\lim_{i\rightarrow \infty}\frac{ d(x_i,y_i)}{d(p, y_i)}=2\sin \frac{d_T(\xi,\eta)}{2}$.
  As $h$ is   Tits induced by  $f$,      
   $\{f(x_i)\}$ Tits converges to $h(\xi)$ and $\{f(y_i)\}$ Tits converges to $h(\eta)$.      Therefore
 $\widetilde{\angle_{f(p)}}(f(x_i),f(y_i))   \rightarrow d_T(h(\xi),h(\eta))$ as $i\rightarrow \infty$.  On the other hand,  
   by considering the comparison triangle of $\triangle f(p)f(x_i)f(y_i)$ we see 
$\sin  \widetilde {\angle_{f(p)}}(f(x_i),f(y_i)) \le \frac{d(f(x_i),f(y_i))}{d(f(p),f(y_i))}$.   It follows that

\vspace{3mm}

   \(\begin{array}{l}
  \sin[d_T(h(\xi), h(\eta))]\\  \\
=\sin[ \lim_{i \rightarrow \infty}\widetilde {\angle_{f(p)}}(f(x_i),f(y_i)) ] 
    =\lim_{i\rightarrow \infty}\sin [\widetilde {\angle_{f(p)}}(f(x_i),f(y_i))] \le  \\   \\
    \le \limsup_{i\rightarrow \infty}  \frac{d(f(x_i),f(y_i))}{d(f(p),f(y_i))}
  \le  \limsup_{i\rightarrow \infty} \frac{L d(x_i,y_i)+A}{\frac{d(p, y_i)}{L} -A} =\\ \\
  = \limsup_{i\rightarrow \infty}L^2 \frac{d(x_i,y_i)}{d(p, y_i)} 
  =\lim_{i\rightarrow \infty}L^2 \frac{ d(x_i,y_i)}{d(p, y_i)}
=2L^2\sin \frac{d_T(\xi,\eta)}{2} \le \\ \\
\le 2L^2  \frac{d_T(\xi,\eta)}{2} 
=L^2 d_T(\xi,\eta)
\le L^2 \epsilon
  =L^2 \frac{\mu}{L^2}
  =\mu.\end{array}\)

\vspace{3mm}

The choice of $\mu$ now implies   $ d_T(h(\xi),h(\eta))\le \lambda \sin [d_T(h(\xi),h(\eta))]\le \lambda   L^2 d_T(\xi,\eta)$.
Similarly we have $d_T(\xi,\eta)\le \lambda   L^2 d_T(h(\xi),h(\eta))$  by considering a quasi-inverse of $f$. 

 \end{proof}

Let $\ol B_i\subset \partial_T X_i$ ($i=1,2$)  be the closure of $B_i$ 
    in $\partial_T X_i$ 
with respect to the Tits metric.  
     The completeness of             $\partial_T X_2$    and  Lemma \ref{l38}
 imply     $h$  uniquely  extends to a        local
bi-Lipschitz      map  from $\ol B_1$   to 
$\partial_T X_2$,  which    we still  denote by  $h$.  
  Since $Core(\partial_T X_2)$ may not be closed in $\partial_T X_2$, 
  $h(\ol B_1 )  $               may not lie in $Core(\partial_T X_2)$.  

\begin{Le}  \label{l40}
{$h: \ol B_1                       \rightarrow \partial_T X_2$ 
  is     Tits  induced by   $f$ and 
 $h(\ol B_1  \cap    Core(\partial_T X_1))   \subset Core(\partial_T X_2)$. }
\end{Le}
\begin{proof}   It  is  easy to see  that $h$     is     Tits  induced by   $f$.  
  Let $\xi\in S  \cap \ol B_1- B_1$    with  
 $S\in {\mathcal C}(\partial_T X_1)$.  %     and $\xi\in S  \cap \ol B_1- B_1$. %We  shall  show $h(\xi)\in g(S)$. 
Fix  some $p\in X_1$ and let $x_i=\gamma_{p\xi}(i)$. Since $\{x_i\}$ Tits converges to $\xi$,
 $\{f(x_i)\}$ Tits converges to $h(\xi)$.  On the other hand, $\{f(x_i)\}\subset f(C_p(S))$  implies
 $L_T(\{f(x_i)\})\subset L_T(f(C_p(S)))=g(S)$. It follows that $h(\xi)\in L_T(\{f(x_i)\})\subset g(S)$.

 \end{proof} 

Lemma  \ref{l40}  implies 
 $h(\ol B_1    \cap    Core(\partial_T X_1) )\subset \ol B_2    \cap    Core(\partial_T X_2) $. 
By considering a quasi-inverse of $f$ we see $h:\ol B_1    \cap    Core(\partial_T X_1)  \rightarrow 
\ol B_2    \cap    Core(\partial_T X_2) $ is a bijective map and its inverse is Tits induced by a 
quasi-inverse of $f$.  
  The proof of  Lemma \ref{l40} 
   shows  $h(S\cap \ol B_1)= g(S)\cap \ol B_2$  for any   
 $S\in {\mathcal C}(\partial_T X_1)$. 
    As  $h_{|S\cap \ol B_1}: S\cap  \ol 
B_1\rightarrow g(S)\cap \ol B_2$
 is  Tits  induced by   $f$,   it  preserves the order of points.

%\begin{Le}\label{l41}
%{  For each    $S\in {\mathcal C}(\partial_T X_1)$,   
 %$h_{|S\cap \bar B_1}: S\cap  \bar 
%B_1\rightarrow g(S)\cap \bar B_2$ preserves the order of points.} 
%\end{Le}

\begin{Le}\label{l42}
{Let    $S\in {\mathcal C}(\partial_T X_1)$   with $S\cap \ol B_1\not=\phi$ 
   and $I$ be  a  component of $S-\ol B_1$. 
  Denote the two endpoints of $I$ by $a, b$. 
  Then $L_T(f(C_x(\bar I)))\subset g(S)$ \e{($x\in X_1$)}  is 
one of the two closed  segments in $g(S)$   with endpoints      $h(a)$ and $h(b)$. Furthermore, the interior of
$L_T(f(C_x(\bar I)))$  is a component of  $g(S)-\ol B_2$. }
\end{Le}
\begin{proof}    Let $I_1$ and $I_2$ be the two components of   $g(S)-\{h(a),  h(b)\}$.  
Since $a, b\in \bar I$  we  have $h(a), h(b)\in L_T(f(C_x(\bar I)))$.   
 By  Lemma  \ref{l33}   $L_T(f(C_x(\bar I)))$ is path  connected.   Hence 
   $L_T(f(C_x(\bar I)))\supset \bar I_1$ or  $L_T(f(C_x(\bar I)))\supset \bar I_2$.
   Without loss of generality  we  assume  $L_T(f(C_x(\bar I)))\supset \bar I_1$.  
  Since   $h^{-1}$ is Tits induced by  a quasi-inverse of 
$f$   and  $\bar I\cap \ol B_1=\{a,b\}$,  
    we conclude $L_T(f(C_x(\bar I)))\cap \ol B_2=\{h(a), h(b)\}$. 
 It  follows that $I_1$ is a component of $g(S)-\ol B_2$.  We shall show 
$L_T(f(C_x(\bar I)))= \bar I_1$.

Suppose  $L_T(f(C_x(\bar I)))\not= \bar I_1$. 
Fix  some $\eta\in  L_T(f(C_x(\bar I)))- \bar I_1$.    Then there is a sequence 
$\{x_i\}\subset C_x(\bar I)$  Tits converging to some $\xi\in I$ so that $\{f(x_i)\}$ Tits converges to $\eta$. 
Let $[\eta, h(a)]$ be the closed subinterval of $\bar I_2$ with  endpoints $\eta$ and $h(a)$.
Similarly define $[\eta, h(b)]\subset \bar I_2$, $[\xi, a], [\xi,b]\subset \bar I$. 
  Notice $L_T(f(C_x([\xi,a])))\supset [\eta, h(a)]$ and $L_T(f(C_x([\xi,b])))\supset [\eta, h(b)]$. 
 Now   $L_T(f(C_x(\bar I)))\supset \bar I_1$ implies  $L_T(f(C_x(\bar I)))\supset g(S)$, a contradiction.

 \end{proof}

\begin{Le} \label{l43}
{ The map $h:\ol B_1\cap   Core(\partial_T X_1)     \rightarrow \ol B_2  \cap   Core(\partial_T X_2)        $ extends
    to a bijective map  $h: Core(\partial_T X_1)\rightarrow 
Core(\partial_T X_2)$.}
\end{Le}
 \begin{proof}   Let $I$ be a component of $ Core(\partial_T X_1)-\ol B_1$. 
If $I=S$ is a circle, then  $S\cap \ol B_1=\phi$. Therefore  $g(S)\cap \ol B_2=\phi$ and $g(S)$ is a 
 component of  $Core(\partial_T X_2) $.  We   define $h_{|S}$ to be a similarity between $S$ and $g(S)$.

Suppose  $I=(a,b)\subset S$ is a  proper subset   of a circle $S$.  
Let us  use the notation in the proof of 
Lemma \ref{l42}. 
In this case we  define $h_{|\bar I}$ to be  the  unique
similarity  from $\bar I$ to $\bar I_1$         extending $h_{|\{a,b\}}$.   The map 
  $h: Core(\partial_T X_1)\rightarrow 
Core(\partial_T X_2)$  is clearly    bijective.

 \end{proof}

\begin{Le}\label{l44}
 {  Let   $I$    be  a  connected  component of $ Core(\partial_T X_1)-\ol B_1$.   Then   the following holds: 
$length(I)/{L^2}\le length(h(I)) \le L^2 length(I)$.}
\end{Le}
\begin{proof}   By considering a quasi-inverse of $f$, it is sufficient   to prove  $ length(h(I)) \le \lambda L^2 length(I)$
   for any $\lambda>1$.    For    any  fixed  $\lambda >1$,   
choose $\mu$ and $\epsilon$ as in the proof of Lemma \ref{l39}.
  We consider  two cases depending on whether $I$ is a circle.

  First suppose  $I=(a,b)\subset S$ is a proper subset of  a circle $S\subset \partial_T X_1$.    
 Let $a=\xi_0$,    $  \xi_1,\cdots, \xi_n=b$ divide $\bar I$ into subintervals of   equal  length
$< \epsilon$.   Fix  a point $x\in X_1$   and let 
   $x_i^j=\gamma_{x\xi_i}(t_j)$ ($0\le i\le n$)  with $t_j\ra \i$. 
After passing to a subsequence we may assume  $\{f(x_i^j)\}_{j=1}^\infty$ Tits 
 converges  to some  $\eta_i\in   L_T(f(C_x(\bar I)))= h(\bar I)\subset g(S)$.  Notice 
$\eta_0=h(a)$, $\eta_n=h(b)$ since $h_{|\ol B_1\cap   Core(\partial_T X_1)}    $ is Tits induced by $f$. 
Since $d_T(\xi_i, \xi_{i+1})\le \epsilon$,  the proof  of Lemma \ref{l39} 
shows   $d_T(\eta_i, \eta_{i+1})\le \lambda L^2 d_T(\xi_i, \xi_{i+1})$.   
 Thus $length(h(I)) \le \sum^{n-1}_{i=0} d_T(\eta_i,\eta_{i+1})\le \sum^{n-1}_{i=0}\lambda L^2 d_T(\xi_i,\xi_{i+1})=\lambda L^2  length(I)$.

Now assume  $I=S$ is a circle.   As above    let  $\xi_0$,   $\xi_1,\cdots, \xi_n=\xi_0$ divide $S$ into subintervals of 
     equal length
$l< \epsilon$.  The same proof   yields  $\eta_i\in g(S)$($\eta_{n}=\eta_0$)
  with    $d_T(\eta_i, \eta_{i+1})\le \lambda L^2 l$.  Fix some $y\in X_2$. 
  For any $\eta\in g(S)$, let 
$y^j=\gamma_{y\eta}(t_j)$   with $t_j\ra \i$.  We may
    assume $\{f^{-1}(y^j)\}$ Tits converges
  to some $\xi\in S$,  where $f^{-1}$ is a quasi-inverse of $f$.  There is some $i$ with $d_T(\xi, \xi_i)\le {l}/{2}$. 
    Then  $d_T(\eta, \eta_i)\le    \lambda L^2 d_T(\xi, \xi_i)\le  \lambda L^2 \times
 {l}/{2}$.   It follows that the $n$ closed intervals   centered at $\eta_i$ ($i=0,1, \cdots, n-1$)  with length $l \lambda  L^2$ cover
the circle $g(S)$. Therefore $length(g(S))\le nl \lambda  L^2= \lambda  L^2 length(S)$. 

\end{proof}

\begin{Le}\label{l45}
{The   bijective  map
  $h: Core(\partial_T X_1)\rightarrow 
Core(\partial_T X_2)$  is   $L^2$-bi-Lipschitz   with respect to the metric $d_c$. }
 \end{Le}
\begin{proof}  It suffices to  show $h$ is  $\lambda L^2$-Lipschitz   with respect to  $d_c$ for all $\lambda>1$. 
Fix  $\lambda>1$.  Choose $\mu$ and $\epsilon$ as in the proof of Lemma \ref{l39}. 
We claim for any $\xi_1$,     $\xi_2\in  Core(\partial_T X_1)$, if $d_c(\xi_1, \xi_2)< \epsilon$ 
then $d_T(h(\xi_1), h(\xi_2))\le  \lambda L^2 d_c(\xi_1, \xi_2)$.
  Let us assume the claim for a while.  For any $\xi, \eta\in Core(\p_T X_1)$   with 
 $d_c(\xi,\eta)<\i$,  let $\alpha: [a, b]\rightarrow  Core(\partial_T X_1)$ 
  be  a minimal Tits geodesic  contained in  $Core(\p_T X_1)$   and 
 connecting   $\xi$ and  $\eta$. 
 The claim implies  the map $h\circ \alpha $ is rectifiable
  with  respect to the Tits metric $d_T$  and 
$length(h\circ \alpha)\le \lambda L^2 length(\alpha)= \lambda L^2 d_c(\xi,\eta)  $, 
  hence  $d_c(h(\xi),h(\eta))\le  \lambda L^2 d_c(\xi,\eta)  $.

%a minimal Tits geodesic  contained in  $Core(\p_T X)$   and 
%connecting   $\xi$ and  $\eta$. 

Now we prove the claim.

\noindent
{{\e{Case 1.}}}  $\xi_1,  \xi_2   \in \ol B_1$.  In this case the claim follows from 
  the definition of  $h$   and the proof of 
Lemma \ref{l39}
 since $d_T(\xi_1, \xi_2)=d_c(\xi_1, \xi_2)$   when   $d_c(\xi_1, \xi_2)< \i$. 

\noindent
{\e{Case 2.}} 
  Exactly   one of $\xi_1$, $\xi_2$ lies in  $\ol B_1$, 
  say $\xi_1\in \ol B_1$ and $\xi_2 \notin  \ol B_1$.
   $d_c(\xi_1, \xi_2)< \epsilon$  implies   $\xi_1\xi_2\subset Core(\partial_T X_1)$.
  Since   $\xi_1\xi_2\cap  \ol B_1$ is closed in $\xi_1\xi_2$    and   $\xi_2 \notin  \ol B_1$,
   there is a point $a\in  \xi_1\xi_2$ with $d_T(a,\xi_2)>0$ and
 $a\xi_2\cap \ol B_1=\{a\}$.   Then there is a component $I$ of $Core(\partial_T X_1)-\ol B_1$
   with  $a\xi_2\subset \ol I$. 
%The 
  % minimal  geodesic segment connecting
%$\xi_1$ and $\xi_2$  must pass through $a$ or $b$ since there is no branch point on $I$.  We may assume it passes through $a$. 
  We also have  $d_c(\xi_1, \xi_2)=d_c(\xi_1, a)+d_c(a, \xi_2)$.  Since the map $h_{|I}$ is a similarity and
 $length(h(I))\le    L^2 length(I)$, we have $d_T(h(a), h(\xi_2))\le L^2 d_T(a, \xi_2)$. On the other hand, 
 $d_T(h(\xi_1), h(a))\le \lambda L^2 d_T(\xi_1, a)$ by Case 1.   Therefore  
$d_T(h(\xi_1), h(\xi_2))\le  d_T(h(\xi_1), h(a))+ d_T(h(a), h(\xi_2))\le \lambda L^2 d_T(\xi_1, a)+  L^2 d_T(a, \xi_2)
\le \lambda L^2 (d_T(\xi_1, a)+ d_T(a, \xi_2))=
 \lambda L^2 (d_c(\xi_1, a)+ d_c(a, \xi_2))= \lambda L^2 d_c(\xi_1, \xi_2)$.  

\noindent
{\e{Case 3.}}  $\xi_1, \xi_2 \notin  \ol B_1$. 
 The proof in this case   is similar to   that  in Case 2. 

\end{proof}

Lemma   \ref{l45}  completes the proof  of Theorem  \ref{t17}.

%\textbf{ Acknowledgment} \textit{This paper contains some  results from the author's 
%dissertation
%at the University of Utah. The author would like to thank his advisor,
 % Bruce Kleiner, for numerous suggestions and discussions, without which this work would never be possible.  In particular, the notion of support set is due to him.} 

%\cite{B, BB, BGS, BH, BS, BW, CK, CK2, G, GH, K, K2, KKL, KL, KO, L, N,  W}

 \addcontentsline{toc}{subsection}{References}
%\noindent 
%Xiangdong Xie\\
%Department of Mathematics\\
%Washington University \\
%St.Louis, MO 63130\\
%xxie@math.wustl.edu\\
\end{document}